\documentclass[11pt,a4paper]{article}
\usepackage{amssymb,amsmath}
\newcommand{\trp}{{\sf\scriptsize T}}
\newcommand\eop{\hfill$\Box$\\}
\newtheorem{theorem}{Theorem}[section]
\newtheorem{lemma}[theorem]{Lemma}
\newtheorem{corollary}[theorem]{Corollary}
\newcommand\asto{\,\stackrel{\rm\scriptsize a.s.}{\longrightarrow}\,}
\newcommand\dto{\,\stackrel{\rm\scriptsize d}{\longrightarrow}\,}
\newcommand\Ifkt {\hskip.1em 1\hskip-.6em 1\hskip.1em}

\begin{document}
\title{The adaptive 
Wynn-algorithm in generalized linear models\\ with univariate response}
\author{Fritjof Freise$^1$, Norbert Gaffke$^2$, \ and \  Rainer Schwabe$^2$ \\[2ex]
\phantom{xxxxxxxxxx}\parbox{11cm}{\small $^1$TU Dortmund University \ and \ 
$^2$University of Magdeburg}
}
\date{}
\maketitle

\begin{abstract}
For a nonlinear regression model the information matrices of designs depend on the parameter
of the model. The adaptive Wynn-algorithm for D-optimal design estimates the parameter
at each step on the basis of the employed design points and observed responses so far, and selects 
the next design point as in the classical Wynn-algorithm for D-optimal design. The name `Wynn-algorithm'
is in honor of Henry P.~Wynn who established the latter `classical' algorithm in his 1970 paper \cite{Wynn}.
The asymptotics of the sequences of designs and maximum likelihood estimates generated by the adaptive algorithm 
is studied for an important class of nonlinear regression models: generalized linear models 
whose (univariate) response variables follow a distribution from a one-parameter exponential family. 
Under the assumptions of compactness of the experimental region and of the parameter space 
together with some natural continuity assumptions it is shown that the  
adaptive ML-estimators are strongly consistent and the design sequence is 
asymptotically locally D-optimal at the true parameter point.
If the true parameter point  is an interior point of the parameter space then
under some smoothness assumptions the asymptotic normality of the adaptive ML-estimators is obtained.
\end{abstract} 

\section{Introduction}
\setcounter{equation}{0}
In a nonlinear regression model the information matrix of a design depends on the model
parameter $\theta\in\Theta$ whose true value is unknown. Modifying the classical algorithm
of Wynn \cite{Wynn} for sequential generation of a D-optimal design in linear regression
to an adaptive sequential procedure in a nonlinear model, the `adaptive Wynn-algorithm' emerges,
which was called  `one-step ahead adaptive D-optimal design algorithm' in Pronzato \cite{Pronzato}.

By $\mathbb{N}$, $\mathbb{N}_0$, $\mathbb{R}$, and $\mathbb{R}^p$  we denote the set of all 
positive integers, the set of all nonnegative integers, the real line, 
and the $p$-dimensional Euclidean space, respectively. 
Vectors $a\in\mathbb{R}^p$ are written as column vectors and $a^\trp$ denotes the transposed of $a$, which is
a $p$-dimensional row vector. The usual Euclidean norm on $\mathbb{R}^p$ is denoted by $\Vert a\Vert=(a^\trp a)^{1/2}$.
If $(a_i)_{i\in I}$ is a family of vectors $a_i\in \mathbb{R}^p$ then ${\rm span}\bigl\{a_i\,:\,i\in I\bigr\}$ denotes 
the linear subspace of $\mathbb{R}^p$ generated by the vectors $a_i$ ($i\in I$). 
For a linear subspace $V$ of $\mathbb{R}^p$ the dimension of $V$ is denoted by ${\rm dim}(V)$.
If $A$ is a symmetric
$p\times p$ matrix then ${\rm tr}(A)$ denotes the trace of $A$ and $\Vert A\Vert$ denotes the Frobenius norm of $A$, i.e., 
$\Vert A\Vert=\bigl({\rm tr}(A^2)\bigr)^{1/2}$. For any two symmetric $p\times p$ matrices $A$ and $B$
we write $A\le B$ or, equivalently, $B\ge A$ iff $B-A$ is nonnegative definite. Thereby a semi-ordering is defined on  
the set of all symmetric $p\times p$ matrices, which is called the Loewner semi-ordering. 

We give an outline of the adaptive Wynn-algorithm.
Let ${\cal X}$ be the experimental region and $\Theta$ be the parameter space.
For each $\theta\in\Theta$ a function \ $f_\theta\,:\,{\cal X}\longrightarrow\mathbb{R}^p$ is given such that
the range of $f_\theta$ spans $\mathbb{R}^p$,
i.e., ${\rm span}\bigl\{f_\theta(x)\,:\,x\in{\cal X}\bigr\}=\mathbb{R}^p$ for each $\theta\in\Theta$.
Throughout it is assumed that ${\cal X}$ and $\Theta$ are compact metric spaces with distance functions
$d_{{\cal X}}$ and $d_\Theta$, resp., and  
the function $(x,\theta)\mapsto f_\theta(x)$ is continuous on ${\cal X}\times \Theta$.
Of course,  the assumption of compactness of the parameter space $\Theta$ is somewhat disturbing but, presently,   
indispensable for our results. However,  
in the literature on adaptive procedures in stochastic approximation it is quite common to assume
compactness of the parameter space and, moreover, to assume the true parameter to be an interior point,
see e.g. Venter \cite{Venter}, Section 4.

An (approximate) design $\xi$ is a probability measure with finite support on ${\cal X}$,
and it can formally be represented as
\[
\xi\,=\,\sum_{x\in{\rm\scriptsize supp}(\xi)}\xi(x)\,\delta_x ,
\]
where ${\rm supp}(\xi)$ denotes the support of $\xi$, which is a nonempty finite subset of ${\cal X}$, and to each $x\in{\rm supp}(\xi)$ 
the design $\xi$ assigns a positive weight $\xi(x)>0$ such that $\sum_{x\in{\scriptsize\rm supp}(\xi)}\xi(x)=1$. 
The symbol $\delta_x$ (for any $x\in{\cal X}$) stands for the one-point probability measure on ${\cal X}$ concentrated at the point $x$.
For a design $\xi$ and for a parameter point $\theta\in\Theta$ the information matrix (per observation) of $\xi$ at $\theta$
is given by
\begin{equation}
M(\xi,\theta)\,=\,\sum_{x\in{\rm\scriptsize supp}(\xi)} \xi(x)\,f_\theta(x)f_\theta^\trp(x)\ ,\label{eq1-0}
\end{equation}
which is a nonnegative definite $p\times p$ matrix. 
The information matrices defined by 
(\ref{eq1-0}) arise as Fisher information in some nonlinear regression model and, in particular, 
the functions $f_\theta$ are related to a local linearization at $\theta$ of the (univariate) 
nonlinear mean response $\mu(x,\theta)$, say.
E.g., in case of a homoscedastic regression model the vector $f_\theta(x)$ is given by the gradient of 
$\mu(x,\,\cdot\,)$ at $\theta$. In case of heteroscedasticity, also the variance function and possibly its gradient
enters into $f_\theta$, see Atkinson et al.~\cite{Atkinson-et-al}.  
For the case of a generalized linear model
the functions $f_\theta$ have the pleasant property that
the parameter $\theta$ only enters into a positive scalar factor, i.e., a real-valued positive function $\psi(x,\theta)$ 
while the `body' of the functions $f_\theta$ is given by one $\mathbb{R}^p$-valued function $f$. We will refer to this
situation as `condition (GLM)` on the family of functions $f_\theta$, $\theta\in\Theta$, namely:\\[.5ex] 
{\bf Condition (GLM)}\\
$f_\theta(x)\,=\,\psi(x,\theta)\,f(x)$ \ for all $(x,\theta)\in {\cal X}\times \Theta$, where
$\psi\,:\,{\cal X}\times\Theta\longrightarrow(\,0\,,\,\infty)$ \ and \ $f:\,{\cal X}\longrightarrow\mathbb{R}^p$ \ 
are given continuous functions. 
\eop
\indent
For a generalized linear model one has, even more specially, that $\Theta\subseteq\mathbb{R}^p$ and
the real-valued function $\psi$ is actually a function of $f^\trp(x)\,\theta$, i.e.,
\begin{equation}
\psi(x,\theta)\,=\,\varphi\bigl(f^\trp(x)\,\theta\bigr),\ \ x\in{\cal X},\ \theta\in\Theta,\label{eq1-0b}  
\end{equation}
where $\varphi$ is a continuous function of one real variable. As an example, for the logistic model with Bernoulli response 
variables one has
\begin{eqnarray*}
&&\mu(x,\theta)\,=\,\exp\bigl(f^\trp(x)\,\theta\bigr)\Big/\Bigl(1+\exp\bigl(f^\trp(x)\,\theta\bigr)\Bigr),\ \
x\in{\cal X},\ \theta\in\Theta,\\
&&\mbox{and }\ 
 \varphi(u)\,=\,\exp(u/2)\,\big/\,\bigl(1 + \exp(u)\bigr),\ \ u\in\mathbb{R}.
\end{eqnarray*}
see Atkinson and Woods \cite{Atkinson-Woods}, Section 2.3.
     
The adaptive Wynn algorithm generates a sequence of designs 
$\xi_n$, $n\ge n_{\rm\scriptsize st}$, (the index `st' standing for `starting')
which is obtained from a sequence of points $x_i\in{\cal X}$, $i\in\mathbb{N}$, and a sequence of parameter points 
$\theta_n\in\Theta$, $n\ge n_{\rm\scriptsize st}$, as follows,
\begin{eqnarray}
&&\xi_n\,=\,\frac{1}{n}\sum_{i=1}^n\delta_{x_i}\quad\mbox{for all $n\ge n_{\rm\scriptsize st}$,}\label{eq1-1}\\
&&x_{n+1}\,=\,\arg\max_{x\in{\cal X}}f_{\theta_n}^\trp(x)M^{-1}(\xi_n,\theta_n)f_{\theta_n}(x)
\quad\mbox{for all $n\ge n_{\rm\scriptsize st}$,} \label{eq1-2}  
\end{eqnarray}
where it is assumed that the starting design $\xi_{n_{\rm\tiny st}}$ is such that its information matrix 
$M(\xi_{n_{\rm\tiny st}},\theta)$ is positive definite for all $\theta\in\Theta$.
This implies
positive definiteness of the information matrices of all designs $\xi_n$, $n\ge n_{\rm\scriptsize st}$, since 
\begin{eqnarray}
\xi_{n+1}&=&{\textstyle\frac{n}{n+1}}\xi_n + {\textstyle\frac{1}{n+1}}\delta_{x_{n+1}}\ ,\quad\mbox{and hence}\label{eq1-3}\\
M\bigl(\xi_{n+1},\theta\bigr)&=&{\textstyle\frac{n}{n+1}}M\bigl(\xi_n,\theta\bigr)\,+\,
{\textstyle\frac{1}{n+1}}f_\theta(x_{n+1})f^\trp_\theta(x_{n+1}),\nonumber
\end{eqnarray} 
which entails \ $M(\xi_n,\theta)\ge(n_{\rm\scriptsize st}/n)\,M(\xi_{n_{\rm\tiny st}},\theta)$, 
for all $n\ge n_{\rm\scriptsize st}$ and all $\theta\in\Theta$. 
Note that  the design
 $\xi_n$ for each  $n\ge n_{\rm\scriptsize st}$ is an exact design of size $n$
since the weights assigned to its support points are integer multiples of $1/n$, and hence $\xi_n$ can be
exactly realized for the sample size $n$.  The sequence of parameter points $\theta_n\in\Theta$, $n\ge n_{\rm\scriptsize st}$, 
employed will actually be generated by adaptive parameter estimation, i.e.,
$\theta_n=\widehat{\theta}_n(x_1,y_1,\ldots,x_n,y_n)$ for all $n\ge n_{\rm\scriptsize st}$,
where $y_1,\ldots y_n,\ldots$ are the sequentially observed univariate responses at the design points 
$x_1,\ldots,x_n,\ldots$, resp., due to an underlying regression model with a mean response 
function $\mu(x,\theta)$, $x\in{\cal X}$, $\theta\in\Theta$, as mentioned above.  

In Section 2 we study the asymptotic behavior of the design sequence $\xi_n$ 
and their information matrices under any sequence of parameter points $\theta_n$,  $n\ge n_{\rm\scriptsize st}$,
which may be thought of as a path of a sequence of adaptive estimators $\widehat{\theta}_n$, 
$n\ge n_{\rm\scriptsize st}$. Also, the design sequence $\xi_n$, $n\ge n_{\rm\scriptsize st}$,
may be viewed as a path of a sequence of adaptive random designs.   
In Section 3 the asymptotic properties (strong consistency, asymptotic normality) of adaptive ML-estimators in the algorithm 
are derived. For modelling the adaptive procedure inherent in the algorithm we follow the martingale approach of 
Lai and Wei \cite{Lai-Wei}, Lai \cite{Lai}, and Chen,  Hu and Ying \cite{Chen-Hu-Ying}.  
Some known results on matrices used in our proofs
are collected in the appendix.
   
The paper of Pronzato  \cite{Pronzato} 
deals with the adaptive Wynn-algorithm for the case of a finite design space (and a compact parameter space).
In that paper, under some conditions of Chebyshev type on the functions $f_\theta$, $\theta\in\Theta$, and the mean response function,
asymptotic results of the design sequence and of adaptive least squares estimators were derived, and also for 
adaptive ML-estimators in the particular case of binary response variables.     
The thesis of Freise \cite{Fritjof} provides an interesting contribution to the asymptotics of the adaptive Wynn algorithm.  
Of further interest, though not dealing with adaptive procedures,
are the papers of Wu \cite{Wu} on nonlinear least squares estimators, and 
of Fahrmeir and Kaufmann \cite{Fahrmeir-Kaufmann} on maximum likelihood estimators
in generalized linear models.

\section{Asymptotic properties of designs}
\setcounter{equation}{0}     
Throughout this section let $\theta_n\in\Theta$,  $n\ge n_{\rm\scriptsize st}$, be any given sequence of parameter points
and let $\xi_n$, $n\ge n_{\rm\scriptsize st}$, be the sequence of designs given by (\ref{eq1-1}) and (\ref{eq1-2}),
where the starting design $\xi_{n_{\rm\tiny st}}$ is such that its information matrix $M(\xi_{n_{\rm\tiny st}},\theta)$
is positive definite for all $\theta\in\Theta$, and hence $M(\xi_n,\theta)$ is positive definite for all
$n\ge n_{\rm\scriptsize st}$ and all $\theta\in\Theta$. 

An important question is whether positive definiteness of the information matrices of the designs
$\xi_n$ is preserved asymptotically in the sense that 
\begin{equation}
\inf_{n\ge n_{\rm\tiny st}} \lambda_{\scriptsize\rm min}\bigl(M(\xi_n,\theta_n)\bigr)\ >0,\label{eq2-1a}
\end{equation}
or, even stronger, 
\begin{equation} 
\inf_{n\ge n_{\rm\tiny st},\,\theta\in\Theta} \lambda_{\scriptsize\rm min}\bigl(M(\xi_n,\theta)\bigr)\ >0, \label{eq2-1b}
\end{equation}  
where $\lambda_{\rm\scriptsize min}(A)$ denotes the smallest eigenvalue of a symmetric matrix $A$.
Answers to the questions about asymptotic nonsingularity will be given. 
Under condition (GLM) the stronger asymptotic nonsingularity (\ref{eq2-1b}) 
holds true, while a weaker technical condition (T) ensures the asymptotic nonsingularity (\ref{eq2-1a}). 
We start our derivations with four lemmas. 

For a real number $a$ we denote by \  $\lceil a\rceil$ \ the smallest integer greater than or equal to $a$.
\vspace*{1ex}
\begin{lemma}\quad\label{lem2-1}\\     
Let $\beta_n$, $n\ge m_0$, be a sequence in $[\,0\,,\,1]$, where $m_0\in\mathbb{N}$ is given,  
and let $\beta\in(\,0\,,\,1)$ such that for each $n\ge m_0$ the following two implications hold.
\begin{eqnarray} 
&&\mbox{If \ $\beta_n>\beta$ \ then \ $\beta_{n+1}\,=\,\frac{n}{n+1}\,\beta_n$\,;}\label{eq2-2a}\\
&&\mbox{if \ $\beta_n\le\beta$ \ then \ $\beta_{n+1}\,\le\,\beta_n+\frac{1}{n+1}$\,.}\label{eq2-3a}
\end{eqnarray}
Let $\widetilde{\beta}>\beta$ be given. Denote \ $m_1\,=\,m_1(\beta,\widetilde{\beta},m_0)\,:=\,
\big\lceil 1/\beta\big\rceil\,\max\bigl\{m_0,\,\big\lceil1/(\widetilde{\beta}-\beta)\big\rceil\bigr\}$.\\ 
Then: \  $\beta_n\,\le\,\widetilde{\beta}$ \ for all $n\ge m_1$.
\end{lemma}

\noindent
{\bf Proof.} \ 
We show that 
\begin{equation}
\mbox{if $n_1\ge m_0$ and $\beta_{n_1}\le\beta$ then \ $\beta_n\le \beta+\frac{1}{n_1}$ for all $n\ge n_1$.}
\label{eqB-7}
\end{equation}   
Let $n_1\ge m_0$ with $\beta_{n_1}\le\beta$ be given. In case that the sequence $\beta_n$, $n\ge n_1$, never exceeds $\beta$ 
the conclusion in (\ref{eqB-7}) trivially holds. 
In the other case, by (\ref{eq2-2a}), it suffices to show that $\beta_n\le \beta+\frac{1}{n_1}$ holds for those $n> n_1$ for which 
$\beta_{n-1}\le \beta$ and $\beta_n>\beta$. For such $n$, by (\ref{eq2-3a}), $\beta_n\le \beta_{n-1}+\frac{1}{n}\le \beta+\frac{1}{n_1}$.\\
Next we show that 
\begin{equation}
\mbox{if $n_2\ge m_0$ and $\beta_{n_2}>\beta$ then $\beta_\nu\le\beta$ for some $\nu\in\{n_2+1,\ldots,\lceil n_2/\beta\rceil\}$.}
\label{eqB-8}
\end{equation}
Let $n_2\ge m_0$ with $\beta_{n_2}>\beta$ be given. If $r$ is a nonnegative  integer such that $\beta_{n_2+k}>\beta$ for all $k=0,\ldots,r$, 
then by (\ref{eq2-2a}) 
\ $\beta<\beta_{n_2+r}=\frac{n_2}{n_2+r}\beta_{n_2}\le\frac{n_2}{n_2+r}$ and hence $\beta<\frac{n_2}{n_2+r}$,
i.e., $n_2+r<n_2/\beta$. So there must be some $\nu\in\{n_2+1,\ldots,\lceil n_2/\beta\rceil\}$ such that $\beta_{\nu}\le\beta$.\\
Consider $m_1$ as defined in the lemma and define 
$k_1=\max\{m_0,\lceil 1/(\widetilde{\beta}-\beta)\rceil\}$. 
Note that $m_1=\big\lceil 1/\beta\big\rceil\,k_1$. We show that $\beta_n\le\widetilde{\beta}$
for all $n\ge m_1$.\\        
\underline{Case 1:} $\beta_{k_1}\le\beta$. By (\ref{eqB-7}) with $n_1=k_1$ one gets 
$\beta_n\le \beta+\frac{1}{k_1}\le \widetilde{\beta}$ for all $n\ge m_1$.\\
\underline{Case 2:} $\beta_{k_1}>\beta$. By (\ref{eqB-8}) with $n_2=k_1$ one gets some 
$\nu\in\{k_1+1,\ldots,\lceil k_1/\beta\rceil\}$ such that $\beta_\nu\le\beta$. Application of (\ref{eqB-7}) on $n_1=\nu$
yields $\beta_n\le \beta+\frac{1}{\nu}\le \widetilde{\beta}$ for all $n\ge \nu$ and, in particular, 
 $\beta_n\le \widetilde{\beta}$ for all $n\ge m_1$ since
$\nu\le\big\lceil k_1/\beta\big\rceil \le \lceil 1/\beta\rceil\,k_1=m_1$. 
\eop

We will use the two positive real constants given by
\begin{eqnarray}
\gamma &:=&\sup_{x\in{\cal X},\, \theta\in\Theta}\Vert f_\theta(x)\Vert\ ,\label{eq2-2}\\
\kappa &:=& \inf_{\Vert v\Vert=1,\, \theta\in\Theta}\ \max_{x\in{\cal X}}\bigl(v^\trp f_\theta(x)\bigr)^2,
\label{eq2-3}
\end{eqnarray}
where in (\ref{eq2-3}) the infimum  is taken over all $v$ from the unit sphere of $\mathbb{R}^p$ and over all $\theta\in\Theta$. 
In fact, both the supremum in (\ref{eq2-2}) and the infimum in (\ref{eq2-3}) are attained and are positive. 
This is obvious for the former supremum by the continuity and compactness assumptions.
For the infimum in  (\ref{eq2-3}), note that 
the function $(v,\theta)\mapsto \max_{x\in{\cal X}}\bigl(v^\trp f_\theta(x)\bigr)^2$ is lower semi-continuous 
(as a pointwise maximum of a family of continuous functions)
and positive, where the latter follows from the basic assumption that the image $\{f_\theta(x)\,:\,x\in{\cal X}\}$ spans 
$\mathbb{R}^p$ for each $\theta\in\Theta$. By
compactness of the unit sphere of $\mathbb{R}^p$ and compactness of $\Theta$ the infimum in (\ref{eq2-3}) 
is attained and hence positive.

\begin{lemma}\quad\label{lem2-2}\\
Let $\xi_0$ be a design and $\theta_0\in\Theta$ such that $M(\xi_0,\theta_0)$ is positive definite.\\
Let $x_0\,=\,\arg\max_{x\in{\cal X}}f_{\theta_0}^\trp(x)M^{-1}(\xi_0,\theta_0)f_{\theta_0}(x)$ and $\eta\in(\,0\,,\,1\,)$.
Then for all $a\in\mathbb{R}^p$ such that \ $\Vert f_{\theta_0}(x_0)-a\Vert\le\eta\kappa/\gamma$ one has
\[
f_{\theta_0}^\trp(x_0)M^{-1}(\xi_0,\theta_0)f_{\theta_0}(x_0)\,\le\,\frac{1}{(1-\eta)^2}\,a^\trp M^{-1}(\xi_0,\theta_0)\,a.
\]
\end{lemma}

\noindent
{\bf Proof.} \ 
Abbreviate $M_0=M(\xi_0,\theta_0)$. Define  \ $b_0\,:=\,M^{-1}_0f_{\theta_0}(x_0)\big/
\bigl(f_{\theta_0}^\trp(x_0)M^{-1}_0f_{\theta_0}(x_0)\bigr)$. \ Then
\begin{equation}
\Vert b_0\Vert^2\,=\,
\frac{f_{\theta_0}^\trp(x_0)M^{-2}_0f_{\theta_0}(x_0)}{
\bigl(f_{\theta_0}^\trp(x_0)M^{-1}_0f_{\theta_0}(x_0)\bigr)^2}. \label{eqB-11a}
\end{equation}
Denote by  $\lambda_1$ the smallest eigenvalue of $M_0$. We will show that
\begin{eqnarray}
  f_{\theta_0}^\trp(x_0)M^{-2}_0f_{\theta_0}(x_0) &\le &\lambda_1^{-2}\,\gamma^2\ \ 
\mbox{and}\label{eqB-11}\\   
  f_{\theta_0}^\trp(x_0)M^{-1}_0f_{\theta_0}(x_0) &\ge & \lambda_1^{-1}\,\kappa\ ,
\label{eqB-12}
\end{eqnarray}
To prove (\ref{eqB-11}), the obvious inequality (in the Loewner semi-ordering) \ $M^{-2}_0\le
\lambda_1^{-2}\,I_p$, where $I_p$ denotes the $p\times p$ unit matrix, yields
\[
 f_{\theta_0}^\trp(x_0)M^{-2}_0f_{\theta_0}(x_0)\,\le 
\lambda_1^{-2}\,f_{\theta_0}^\trp(x_0)f_{\theta_0}(x_0)\,=\,
\lambda_1^{-2}\,\Vert f_{\theta_0}(x_0)\Vert^2\le\,\lambda_1^{-2}\,\gamma^2.
\]
To prove (\ref{eqB-12}) let $v_1$ be a normalized eigenvector to $\lambda_1$ of $M_0$.
By definition of $\kappa$ from (\ref{eq2-3})
\[
\kappa\,\le \max_{x\in{\cal X}}\bigl\{\bigl(v_1^\trp f_{\theta_0}(x)\bigr)^2\bigr\}\,=\,
\bigl(v_1^\trp f_{\theta_0}(z_0)\bigr)^2
\] 
for some $z_0\in{\cal X}$. Hence, together with the obvious inequality (in the Loewner semi-ordering) \ 
$M^{-1}_0\ge \lambda_1^{-1}\,v_1v_1^\trp$, one obtains
\begin{eqnarray*}
&&f_{\theta_0}^\trp(x_0)M^{-1}_0f_{\theta_0}(x_0)\,\ge\, 
f_{\theta_0}^\trp(z_0)M^{-1}_0f_{\theta_0}(z_0)\,\ge\,
\lambda_1^{-1}\,f_{\theta_0}^\trp(z_0)v_1v_1^\trp f_{\theta_0}(z_0)\\
&&\,=\,\lambda_1^{-1}\,\bigl(v_1^\trp f_{\theta_0}(z_0)\bigr)^2\ge\,
\lambda_1^{-1}\,\kappa\ .
\end{eqnarray*}
From (\ref{eqB-11a}), (\ref{eqB-11}), and (\ref{eqB-12}) we get
\begin{equation}
\Vert b_0\Vert\ \le \gamma/\kappa . \label{eqB-12a}
\end{equation}
Let $a\in\mathbb{R}^p$ such that 
$\Vert f_{\theta_0}(x_0)-a\Vert\le\eta\kappa/\gamma$. Recall that, by definition of $b_0$,
\[
\frac{1}{b_0^\trp M_0\,b_0}\,=\, f_{\theta_0}^\trp(x_0)\,M^{-1}_0\,f_{\theta_0}(x_0)
\ \mbox{ and }\ f_{\theta_0}^\trp(x_0)\,b_0=1.
\]
Together with (\ref{eqB-12a}) we get
\[
|a^\trp b_0-1|=|a^\trp b_0-f_{\theta_0}^\trp(x_0)\,b_0|=\big|\bigl(a-f_{\theta_0}(x_0)\bigr)^\trp b_0\big|
\le \Vert a-f_{\theta_0}(x_0)\Vert\,\Vert b_0\Vert
\,\le\,\eta,
\]
hence $a^\trp b_0\ge1-\eta>0$. Define $b:=b_0/(a^\trp b_0)$. Then $a^\trp b=1$ and hence by (M3) of the appendix,
\begin{eqnarray*}
&&a^\trp M^{-1}_0\,a\,\ge\,\frac{1}{b^\trp M_0\,b}\,=\,\frac{(a^\trp b_0)^2}{b_0^\trp M_0\,b_0}\\
&&\,=\,(a^\trp b_0)^2\,f_{\theta_0}^\trp(x_0)\,M^{-1}_0\,f_{\theta_0}(x_0)\,\ge\,(1-\eta)^2\,
f_{\theta_0}^\trp(x_0)\,M^{-1}_0\,f_{\theta_0}(x_0),
\end{eqnarray*}
from which the result follows.
\eop
 
\begin{lemma}\quad\label{lem2-3}\\
Let  $p\ge2$ and $\eta\in\bigl(\,0\,,\,1-\frac{1}{\sqrt{p}}\,\bigr)$. Let    
$S\subseteq{\cal X}$ and $n\ge n_{\rm\scriptsize st}$ be given such that 
\[
\Vert f_{\theta_n}(x)-f_{\theta_n}(z)\Vert\,\le\,\eta\kappa/\gamma\ \ \mbox{for all $x,z\in S$, and }\   
\xi_n(S)>\frac{1}{(1-\eta)^2p}. 
\]
Then: \ $x_{n+1}\not\in S$.
\end{lemma}

\noindent{\bf Proof.} \  
Suppose that $x_{n+1}\in S$. Consider the mean (of $f_{\theta_n}(x)$ over $S$ w.r.t. $\xi_n$),
\begin{equation}
\overline{f}_{\theta_n}(S,\xi_n)\,:=\,\frac{1}{\xi_n(S)}\,\sum_{x\in S\cap\,{\rm\scriptsize supp}(\xi_n)}\xi_n(x)\,f_{\theta_n}(x).
\label{eq2-add}
\end{equation}
Since $\Vert f_{\theta_n}(x_{n+1})-f_{\theta_n}(x)\Vert \le \eta\kappa/\gamma$ for all $x\in S$ we get
$\Vert f_{\theta_n}(x_{n+1})-\overline{f}_{\theta_n}(S,\xi_n)\Vert\le\eta\kappa/\gamma$. By Lemma \ref{lem2-2},
\[
f_{\theta_n}^\trp (x_{n+1})\,M^{-1}(\xi_n,\theta_n)\,f_{\theta_n}(x_{n+1})\,\le\,\frac{1}{(1-\eta)^2}\,
\overline{f}_{\theta_n}^\trp(S,\xi_n)\,M^{-1}(\xi_n,\theta_n)\,\overline{f}_{\theta_n}(S,\xi_n).
\]
By (M1) and (M2) of the appendix,
\begin{eqnarray*}
&&\overline{f}_{\theta_n}^\trp(S,\xi_n)\,M^{-1}(\xi_n,\theta_n)\,\overline{f}_{\theta_n}(S,\xi_n)\\
&&\le\
\overline{f}_{\theta_n}^\trp(S,\xi_n)\,
\Bigl(\xi_n(S)\sum_{x\in S\cap\,{\rm\scriptsize supp}(\xi_n)}
{\textstyle\frac{\xi_n(x)}{\xi_n(S)}}
f_{\theta_n}(x)\,f_{\theta_n}^\trp(x)\Bigr)^-\,\overline{f}_{\theta_n}(S,\xi_n)\\
&&\le\ \frac{1}{\xi_n(S)}\,<\,(1-\eta)^2\,p.
\end{eqnarray*}
Hence it follows that 
$f_{\theta_n}^\trp (x_{n+1})\,M^{-1}(\xi_n,\theta_n)\,f_{\theta_n}(x_{n+1})\,<\,p$.
This is a contradiction since we know from the Kiefer-Wolfowitz equivalence Theorem that
$\max_{x\in{\cal X}}f_{\theta_n}^\trp (x)\,M^{-1}(\xi_n,\theta_n)\,f_{\theta_n}(x)\,\ge\,p$.
So $x_{n+1}\not\in S$ must be true.
\eop 

For $a\in\mathbb{R}^p$, $\emptyset\not=V\subseteq\mathbb{R}^p$, and $\varepsilon>0$ we denote 
\[
{\rm dist}(a,V)=\inf_{v\in V}\Vert a-v\Vert \ \mbox{ and }\ 
\overline{V}(\varepsilon)\,:=\,\bigl\{a\in\mathbb{R}^p\,:\,{\rm dist}(a,V)\le \varepsilon\bigr\}.
\] 
The set $\overline{V}(\varepsilon)$ may be called an $\varepsilon$-neighborhood of $V$. 
For any subset $C\subseteq\mathbb{R}^p$ and $\theta\in\Theta$ we denote, as usual,
$f_\theta^{-1}(C)\,=\,\bigl\{x\in{\cal X}\,:\,f_\theta(x)\in C\bigr\}$. 

\begin{lemma}\quad\label{lem2-4}\\
Let $V\subseteq\mathbb{R}^p$ be a linear subspace with ${\rm dim}(V)\le p-1$, and let 
$\delta$ with $0<\delta\le\sqrt{\kappa}$ and $n\ge n_{\rm\scriptsize st}$ be given. 
Then, denoting \ $w_n:=\xi_n\Bigl(f_{\theta_n}^{-1}\bigl(\overline{V}(\delta)\bigr)\Bigr)$, one has
\[
f_{\theta_n}^\trp(x_{n+1})\,M^{-1}(\xi_n,\theta_n)\,f_{\theta_n}(x_{n+1})\,\ge\,\Bigl(1-\bigl(1-\delta^2/\kappa\bigr)w_n\Bigr)^{-1}.
\]
\end{lemma}

\noindent{\bf Proof.} \ For all $x\in{\cal X}$ decompose
\[
f_{\theta_n}(x)\,=\,u(x)\,+\,v(x),\ \mbox{where }\ u(x)\in V\ \mbox{ and }\ v(x)\in V^\perp,
\]
where $V^\perp$ denotes the orthogonal complement of $V$ in $\mathbb{R}^p$. Choose \ 
$x^*\,=\,\arg\max_{x\in{\cal X}}\Vert v(x)\Vert$. \ 
Clearly, for all  $x\in f_{\theta_n}^{-1}\bigl(\overline{V}(\delta)\bigr)$ one has 
$\Vert v(x)\Vert={\rm dist}\bigl(f_{\theta_n}(x),V\bigr)\le\delta$.
On the other hand,
$\Vert v(x^*)\Vert\ge\sqrt{\kappa}$ which can be seen as follows. Since
${\rm dim}(V)\le p-1$ there is some $(p-1)$-dimensional linear subspace $W\subseteq\mathbb{R}^p$ such that $V\subseteq W$.
There is a representation $W=\{a\in\mathbb{R}^p\,:\,c^\trp a=0\}$ for some $c\in\mathbb{R}^p$ with $\Vert c\Vert=1$. 
By ${\rm dist}(a,V)\ge{\rm dist}(a,W)$ for all $a\in\mathbb{R}^p$, and
by definition of $\kappa$ in (\ref{eq2-3}) one gets
\[
\Vert v(x^*)\Vert\,=\,\max_{x\in{\cal X}}{\rm dist}(f_{\theta_n}(x),V)
\,\ge\,\max_{x\in{\cal X}}{\rm dist}(f_{\theta_n}(x),W)\,=\,\max_{x\in{\cal X}}|c^\trp f_{\theta_n}(x)|\,\ge\sqrt{\kappa}.
\] 
Define $b\,:=\,v(x^*)/\Vert v(x^*)\Vert^2$. Clearly, $b^\trp f_{\theta_n}(x^*)=1$ hence by (M3) of the appendix
\begin{equation}
f_{\theta_n}^\trp(x^*)\,M^{-1}(\xi_n,\theta_n)\,f_{\theta_n}(x^*)\,\ge\,\frac{1}{b^\trp M(\xi_n,\theta_n)\,b}.
\label{eq2-add1}
\end{equation} 
Now,\quad  $\displaystyle b^\trp M(\xi_n,\theta_n)\,b\,=\,\sum_{x\in{\rm\scriptsize supp}(\xi_n)}\xi_n(x)\,\bigl(b^\trp f_{\theta_n}(x)\bigr)^2\,=\,
\sum_{x\in{\rm\scriptsize supp}(\xi_n)}\xi_n(x)\,\frac{\bigl(v^\trp(x^*)\,v(x)\bigr)^2}{\Vert v(x^*)\Vert^4}$,\\
and $\bigl(v^\trp(x^*)\,v(x)\bigr)^2\big/\Vert v(x^*)\Vert^4\le \Vert v(x^*)\Vert^2\,\Vert v(x)\Vert^2\big/\Vert v(x^*)\Vert^4
=\Vert v(x)\Vert^2/\Vert v(x^*)\Vert^2\,\le 1$ for all $x\in{\cal X}$. If $x\in f_{\theta_n}^{-1}\bigl(\overline{V}(\delta)\bigr)$
then  \ $\Vert v(x)\Vert^2/\Vert v(x^*)\Vert^2\le \delta^2/\kappa$. Hence,
partitioning ${\rm supp}(\xi_n)$ into ${\rm supp}(\xi_n)\cap f_{\theta_n}^{-1}\bigl(\overline{V}(\delta)\bigr)$ and
${\rm supp}(\xi_n)\setminus f_{\theta_n}^{-1}\bigl(\overline{V}(\delta)\bigr)$, one gets
\[
b^\trp M(\xi_n,\theta_n)\,b \,\le\,
\bigl(\delta^2/\kappa\bigr)\,w_n\,+\,1-w_n\,=\,
1-\bigl(1-\delta^2/\kappa\bigr)\,w_n,
\]
and together with (\ref{eq2-add1}) the result follows.
\eop

We introduce a technical condition (T) which is weaker than (GLM).
It is motivated by the result of Lemma \ref{lem2-5} below.\\[.5ex]

\noindent
{\bf Condition (T)}\\
For each $\delta>0$ there exist an integer $m_0(\delta)\ge n_{\rm\scriptsize st}$ and a $\delta'>0$ such that 
for all $k,\ell\ge m_0(\delta)$ and all linear subspaces $V\subseteq\mathbb{R}^p$
one has \ 
$f_{\theta_k}^{-1}\bigl(\overline{V}(\delta')\bigr)\subseteq f_{\theta_\ell}^{-1}\bigl(\overline{V}(\delta)\bigr)$.
\eop

\begin{lemma}\quad\label{lem2-5}\\
(i) Condition (GLM) implies condition (T).\\
(ii) If \ $\lim_{n\to\infty}\theta_n=\overline{\theta}$ \ for some $\overline{\theta}\in\Theta$ then condition (T) holds.
\end{lemma}

\noindent{\bf Proof.} \ \underbar{Ad (i)}. Assume (GLM). Denote 
\begin{equation}
\psi_{\rm\scriptsize min}=\inf_{(x,\theta)\in{\cal X}\times\Theta}\psi(x,\theta)\quad\mbox{and}\quad
\psi_{\rm\scriptsize max}=\sup_{(x,\theta)\in{\cal X}\times\Theta}\psi(x,\theta).\label{eq2-add10}
\end{equation} 
By compactness and continuity the infimum and the supremum are attained,
and hence $0<\psi_{\rm\scriptsize min}\le\psi_{\rm\scriptsize max}<\infty$. For a given $\delta>0$ choose 
$m_0(\delta)=n_{\rm\scriptsize st}$ and $\delta'=\delta\, \psi_{\rm\scriptsize min}/\psi_{\rm\scriptsize max}$. Let $k,\ell\ge n_{\rm\scriptsize st}$
and a linear subspace $V\subseteq\mathbb{R}^p$ be given. For any $\theta\in\Theta$ and any $\varepsilon >0$ one has
$f_\theta^{-1}\bigl(\overline{V}(\varepsilon)\bigr)=\bigl\{x\in{\cal X}\,:\,{\rm dist}(\psi(x,\theta)\,f(x),V)\le\varepsilon\bigr\}$,
and \ ${\rm dist}\bigl(\psi(x,\theta)\,f(x),V\bigr)=\psi(x,\theta)\,{\rm dist}(f(x),V)$, hence 
\begin{equation}
f_\theta^{-1}\bigl(\overline{V}(\varepsilon)\bigr)=\bigl\{x\in{\cal X}\,:\,{\rm dist}(f(x),V)\le\varepsilon/\psi(x,\theta)\bigr\}.   
\label{eq2-add2}
\end{equation}
For $\theta=\theta_k$ and $\varepsilon=\delta'$ (\ref{eq2-add2}) yields, observing 
$\delta'/\psi(x,\theta_k)\le \delta'/\psi_{\rm\scriptsize min}=\delta/\psi_{\rm\scriptsize max}$, 
\begin{equation}
f_{\theta_k}^{-1}\bigl(\overline{V}(\delta')\bigr)\subseteq\bigl\{x\in{\cal X}\,:\,{\rm dist}(f(x),V)\le\delta/\psi_{\rm\scriptsize max}\bigr\}.
\label{eq2-add3}
\end{equation}
For $\theta=\theta_\ell$ and $\varepsilon=\delta$ (\ref{eq2-add2}) yields, observing 
$\delta/\psi(x,\theta_\ell)\ge \delta/\psi_{\rm\scriptsize max}$,
\begin{equation}
f_{\theta_\ell}^{-1}\bigl(\overline{V}(\delta)\bigr)\supseteq\bigl\{x\in{\cal X}\,:\,{\rm dist}(f(x),V)\le\delta/\psi_{\rm\scriptsize max}\bigr\}.
\label{eq2-add4}
\end{equation} 
From (\ref{eq2-add3}) and (\ref{eq2-add4}) the inclusion 
$f_{\theta_k}^{-1}\bigl(\overline{V}(\delta')\bigr)\subseteq f_{\theta_\ell}^{-1}\bigl(\overline{V}(\delta)\bigr)$ follows.\\[.5ex]
\underbar{Ad (ii)}. Assume that $\lim_{n\to\infty}\theta_n=\overline{\theta}$ for some $\overline{\theta}\in\Theta$.
By compactness of ${\cal X}\times\Theta$ and continuity (hence uniform continuity) of the function $(x,\theta)\longmapsto f_\theta(x)$
the sequence of functions $f_{\theta_n}$, $n\ge n_{\rm\scriptsize st}$, converges to $f_{\overline{\theta}}$ uniformly on ${\cal X}$. 
So, for any given $\delta>0$ there is an $m_0(\delta)\ge n_{\rm\scriptsize st}$ such that 
\begin{equation}
\Vert f_{\theta_k}(x)-f_{\theta_\ell}(x)\Vert\,\le\,\delta/2\quad\mbox{for all $k,\ell\ge m_0(\delta)$ and all $x\in{\cal X}$.}
\label{eq2-add4a}
\end{equation}
Choose $\delta'=\delta/2$. Let $k,\ell\ge m_0(\delta)$ and a linear subspace $V\subseteq\mathbb{R}^p$ be given. 
Using the well-known inequality
\[
\big|{\rm dist}(a,V)-{\rm dist}(b,V)\big|\,\le\,\Vert a-b\Vert\quad\mbox{for all $a,b\in\mathbb{R}^p$},    
\]
one gets from (\ref{eq2-add4a}) that
\begin{equation}
{\rm dist}(f_{\theta_\ell}(x),V)\,\le\,{\rm dist}(f_{\theta_k}(x),V) \,+\,\delta/2\quad\mbox{for all $x\in{\cal X}$.}
\label{eq2-add5}
\end{equation}
From (\ref{eq2-add5}), using $\delta'=\delta/2$, one gets 
$f_{\theta_k}^{-1}\bigl(\overline{V}(\delta')\bigr)\subseteq  f_{\theta_\ell}^{-1}\bigl(\overline{V}(\delta)\bigr)$.
\eop

\begin{theorem}\quad\label{theo2-1}\\
Assume condition (T).  
Then there exist  $n_0\ge n_{\rm\scriptsize st}$, $\varepsilon>0$,  and  $\alpha\in(\,0\,,\,1\,)$  
such that for all $n\ge n_0$ and all $(p-1)$-dimensional linear subspaces $V_{p-1}$ of $\mathbb{R}^p$ 
one has \ $\xi_n\Bigl(f_{\theta_n}^{-1}\bigl(\overline{V}_{p-1}(\varepsilon)\bigr)\Bigr)\le\alpha$.
\end{theorem}

\noindent
{\bf Proof.} \ 
Firstly, consider the (nearly) trivial case $p=1$. 
The only $0$-dimensional linear subspace of $\mathbb{R}^1$ is $\{0\}$, hence 
$\overline{V}_0(\varepsilon)=[-\varepsilon\,,\,\varepsilon]$ for any $\varepsilon>0$. 
From (\ref{eq1-2}) $\bigl(f_{\theta_n}(x_{n+1})\bigr)^2=\max_{x\in{\cal X}}\bigl(f_{\theta_n}(x)\bigr)^2$ 
and by (\ref{eq2-3}) $\bigl(f_{\theta_n}(x_{n+1})\bigr)^2\ge\kappa$ for all $n\ge n_{\rm\scriptsize st}$. 
Choose a $\delta\in(\,0\,,\,\sqrt{\kappa}\,)$ and choose $m_0(\delta)\ge n_{\rm\scriptsize st}$ and $\delta'>0$
according to condition (T). Then $x_{n+1}\not\in f_{\theta_n}^{-1}\bigl([\,-\delta\,,\,\delta\,]\bigr)$ 
for all $n\ge n_{\rm\scriptsize st}$, hence 
\begin{equation}
x_i\not\in\bigcap_{\ell\ge m_0(\delta)} f_{\theta_\ell}^{-1}\bigl([\,-\delta\,,\,\delta\,]\bigr)\quad
\mbox{for all $i\ge m_0(\delta)+1$.} \label{eq2-add6}
\end{equation}
By (T), for all $n\ge m_0(\delta)$ the set $f_{\theta_n}^{-1}\bigl([\,-\delta'\,,\,\delta'\,]\bigr)$
is a subset of the intersection from (\ref{eq2-add6}) and hence 
$x_i\not\in f_{\theta_n}^{-1}\bigl([\,-\delta'\,,\,\delta'\,]\bigr)$ for all $i\ge m_0(\delta)+1$
and all $n\ge m_0(\delta)$. It follows that
\[
\xi_n\Bigl(f_{\theta_n}^{-1}\bigl([\,-\delta'\,,\,\delta'\,]\bigr)\Bigr)\,\le\,\frac{m_0(\delta)}{n}
\quad\mbox{for all $n\ge m_0(\delta)$.}
\]
So, choosing  $n_0=2m_0(\delta)$, $\varepsilon=\delta'$, and $\alpha=1/2$, the statement of the theorem holds
in case $p=1$. In what follows we assume $p\ge2$.
We will prove by induction the following statement ${\rm S}(r)$ for all $r=0,1,\ldots,p-1$.
\begin{itemize}
\item[${\rm S}(r)$]
There exist $\widetilde{n}_r\ge n_{\rm\scriptsize st}$, $\varepsilon_r>0$, and $\alpha_r\in (\,0\,,\,1\,)$ 
such that \ $\xi_n\Bigl(f_{\theta_n}^{-1}\bigl(\overline{V}_r(\varepsilon_r)\bigr)\Bigr)\le\alpha_r$ 
for all $n\ge\widetilde{n}_r$ and all $r$-dimensional linear subspaces $V_r$ of $\mathbb{R}^p$.
\end{itemize}
Then the result will follow by taking $n_0=\widetilde{n}_{p-1}$, $\varepsilon=\varepsilon_{p-1}$,  and $\alpha=\alpha_{p-1}$.
\\[.5ex]   
\underline{$r=0$}. \ The only $0$-dimensional linear subspace of $\mathbb{R}^p$ is the nullspace $V_0=\{0\}$,
and for any $\varepsilon>0$ one has $\overline{V}_0(\varepsilon)=\{a\in\mathbb{R}^p\,:\,\Vert a\Vert \le\varepsilon\}$,
the closed ball centered at zero with radius $\varepsilon$. Choose any $\eta\in\bigl(\,0\,,\,1-\frac{1}{\sqrt{p}}\,\bigr)$
and let $\delta:=\eta\kappa/(2\gamma)$. Choose $m_0(\delta)\ge n_{\rm\scriptsize st}$ and $\delta'>0$ according to condition (T),
and define
\[
S\,:=\,\bigcap_{\ell\ge m_0(\delta)}f_{\theta_\ell}^{-1}\bigl(\overline{V}_0(\delta)\bigr).
\]
Clearly, if $x,z\in f_{\theta_\ell}^{-1}\bigl(\overline{V}_0(\delta)\bigr)$, i.e., $\Vert f_{\theta_\ell}(x)\Vert\le\delta$ and
$\Vert f_{\theta_\ell}(z)\Vert\le\delta$, then $\Vert f_{\theta_\ell}(x)-f_{\theta_\ell}(z)\Vert\le2\delta=\eta\kappa/\gamma$.
So the subset $S$ has the property that if $n\ge m_0(\delta)$ and $x,z\in S$ then $\Vert f_{\theta_n}(x)-f_{\theta_n}(z)\Vert\le\eta\kappa/\gamma$.
By Lemma \ref{lem2-3}, if $n\ge m_0(\delta)$ and $\xi_n(S)>1\big/\bigl((1-\eta)^2p\bigr)$ then $x_{n+1}\not\in S$.  
Choose an $\alpha_0$ with $1\big/\bigl((1-\eta)^2p\bigr)<\alpha_0<1$. The sequence $\beta_n=\xi_n(S)$, $n\ge m_0(\delta)$, 
along with $\beta=1\big/\bigl((1-\eta)^2p\bigr)$ and $\widetilde{\beta}=\alpha_0$, satisfy the assumptions of
Lemma \ref{lem2-1}, and hence by that lemma $\xi_n(S)\le\widetilde{\beta}=\alpha_0$ for all $n\ge 
m_1=m_1\bigl(\beta,\widetilde{\beta},m_0(\delta)\bigr)$. By (T), $f_{\theta_n}^{-1}\bigl(\overline{V}_0(\delta')\bigr)\subseteq S$
for all $n\ge m_0(\delta)$ and hence $\xi_n\Bigl(f_{\theta_n}^{-1}\bigl(\overline{V}_0(\delta')\bigr)\Bigr)\le\alpha_0$ for all
$n\ge m_1$. So statement ${\rm S}(0)$ holds with $\widetilde{n}_0=m_1$, $\varepsilon=\delta'$, and $\alpha_0$ as 
already introduced.\\[.5ex]
\underline{Induction step.} \ Suppose that for some $r\in\{1,\ldots,p-1\}$ statement ${\rm S}(r-1)$ is true,
and let $\widetilde{n}_{r-1}$, $\varepsilon_{r-1}$, and $\alpha_{r-1}$ be chosen as in statement ${\rm S}(r-1)$.      
Since every linear subspace $V_t\subseteq\mathbb{R}^p$ of dimension $t\le r-1$ can be enlarged to
an $(r-1)$-dimensional linear subspace $V_{r-1}\subseteq\mathbb{R}^p$, where $V_t\subseteq V_{r-1}$ and hence 
$\overline{V}_t(\varepsilon_{r-1})\subseteq\overline{V}_{r-1}(\varepsilon_{r-1})$, the assumed statement ${\rm S}(r-1)$
implies the following.
\begin{eqnarray}
&&\mbox{For all $t$-dimensional linear subspaces $V_t\subseteq\mathbb{R}^p$ with $t\le r-1$ and for all $n\ge\widetilde{n}_{r-1}$} \nonumber\\
&&\mbox{one has }\ \ \xi_n\Bigl(f_{\theta_n}^{-1}\bigl(\overline{V}_t(\varepsilon_{r-1})\bigr)\Bigr)\,\le\,\alpha_{r-1}.
\label{eqC-13}
\end{eqnarray}
The rest of the proof of the induction step is lengthy; it is structured into three steps.\\
\underline{Step 1.} \ 
We introduce some sets and constants. 
\begin{eqnarray}
&&{\cal A}_r:=\Bigl\{A=[a_1,\ldots,a_r]\in\mathbb{R}^{p\times r}\,:\,a_j\in\mathbb{R}^p,\  
 \Vert a_j\Vert\le\gamma, \ 1\le j\le r\Bigr\}, \label{eqC-s1-1}\\
&&{\cal A}_r^* := \Bigl\{A\in{\cal A}_r\,:\,\det(A^\trp A)\ge \textstyle\frac{1}{2}\Bigl(\frac{\varepsilon_{r-1}^2}{4}\Bigr)^r\Bigr\}.
\label{eqC-s1-2}
\end{eqnarray}
Obviously, ${\cal A}_r$ and ${\cal A}_r^*$ are compact sets of $p\times r$ matrices.
It is not quite obvious that ${\cal A}_r^*$ is nonempty which can be seen as follows. (\ref{eqC-13}) implies in particular 
that, choosing any $n\ge \widetilde{n}_{r-1}$, the set 
${\cal X}\setminus f_{\theta_n}^{-1}\bigl(\overline{V}_0(\varepsilon_{r-1})\bigr)$ is nonempty, i.e., there is a $z\in{\cal X}$
such that $\Vert f_{\theta_n}(z)\Vert >\varepsilon_{r-1}$. By $\Vert f_{\theta_n}(z)\Vert\le\gamma$ one has $\varepsilon_{r-1}<\gamma$.
Choosing pairwise orthogonal vectors $a_1,\ldots,a_r\in\mathbb{R}^p$ with $\Vert a_j\Vert = \gamma$,
$1\le j\le r$, one gets a matrix $A=[a_1,\ldots,a_r]\in{\cal A}_r$ with $\det(A^\trp A)=\gamma^{2r}>\varepsilon_{r-1}^{2r}$,
hence $A\in{\cal A}_r^*$. Next, denote by 
$\Vert w\Vert_1\,=\,\sum_{j=1}^r|w_j|$ the $\ell^1$-norm  of a vector $w=(w_1,\ldots,w_r)^\trp\in\mathbb{R}^r$
and define
\begin{equation}
c_r\,:=\,\sup\Bigl\{\big\Vert (A^\trp A)^{-1}A^\trp b\big\Vert_1\,:\,A\in{\cal A}_r^*,\ b\in\mathbb{R}^p,\ \Vert b\Vert\le\gamma\Bigr\}.
\label{eqC-s1-3}
\end{equation}
Choosing any $A=[a_1,\ldots,a_r]\in{\cal A}_r^*$ and $b=a_1$ gives $(A^\trp A)^{-1}A^\trp b=e_1=(1,0,\ldots,0)^\trp$,
and $c_r\ge1$ follows. Together with compactness and continuity one has $1\le c_r<\infty$. 
Again by compactness and continuity one can choose a positive integer $K_r$ 
and nonempty subsets $R_1,\ldots,R_{K_r}\subseteq{\cal X}$ such that
\begin{equation}   
{\cal X}=\bigcup_{k=1}^{K_r}R_k\  
\mbox{ and }\ \Vert f_\theta(x)-f_\theta(z)\Vert\,\le\varepsilon_{r-1}/2\ \ \forall\ x,z\in R_k,\ \ 1\le k\le K_r,\ \forall\ \theta\in\Theta.
\label{eqC-s1-4}
\end{equation}
\begin{equation}
\mbox{Choose $\overline{\alpha}_r$ such that }
\ \frac{K_rc_r^2+\alpha_{r-1}}{K_rc_r^2 +1}\,<\,\overline{\alpha}_r\,<1.
\label{eqC-s1-5}
\end{equation}
Note that  $\alpha_{r-1}\,<\,(K_rc_r^2 +\alpha_{r-1})\big/(K_rc_r^2 +1)$, hence  $\alpha_{r-1}<\overline{\alpha}_r$.
Finally, choose  a $\delta>0$ which satisfies the following three conditions,
\begin{eqnarray}
&& 0\,<\,\delta\,<\,\frac{\kappa}{(c_r+1)\gamma} \label{eqC-s1-6}\\
&& \bigl(K_rc_r^2\bigr)^{-1}\,\Bigl(1-(c_r+1){\textstyle\frac{\gamma}{\kappa}}\delta\Bigr)^2(\overline{\alpha}_r-\alpha_{r-1})
\,>\,1-\Bigl(1-{\textstyle\frac{1}{\kappa}}\delta^2\Bigr)\,\overline{\alpha}_r\,,\label{eqC-s1-7}\\
&& \big|\det(A^\trp A)-\det(B^\trp B)\big|\le{\textstyle\frac{1}{2}\Bigl(\frac{\varepsilon_{r-1}^2}{4}\Bigr)^r}\label{eqC-s1-8}\\
&&\mbox{for all }\ A=[a_1,\ldots,a_r]\in{\cal A}_r,\ B=[b_1,\ldots,b_r]\in{\cal A}_r\ \mbox{with }  
\Vert a_j-b_j\Vert\le\delta,\  1\le j\le r.\nonumber
\end{eqnarray}
In fact, such a $\delta$ exists since, firstly, both sides of the inequality (\ref{eqC-s1-7}) are continuous functions of
a real variable $\delta$ and the (strict) inequality (\ref{eqC-s1-7}) holds for $\delta=0$ by (\ref{eqC-s1-5}).
Secondly, (\ref{eqC-s1-8}) is achieved by the uniform continuity
of the function $A\longmapsto\det(A^\trp A)$ on the compact set ${\cal A}_r$ from (\ref{eqC-s1-1}).\\[.5ex]  
\underline{Step 2.} \ With $\delta$ and $\overline{\alpha}_r$ from Step 1 
we show the following:
\begin{itemize}
\item[]
{\em If $V_r\subseteq\mathbb{R}^p$ is an $r$-dimensional linear subspace and $n\ge \widetilde{n}_{r-1}$
such that\\  
$\xi_n\Bigl(f_{\theta_n}^{-1}\bigl(\overline{V}_r(\delta)\bigr)\Bigr)\,>\,\overline{\alpha}_r$, \ then   
\ $x_{n+1}\not\in f_{\theta_n}^{-1}\bigl(\overline{V}_r(\delta)\bigr)$.}
\end{itemize}
Let  an $r$-dimensional linear subspace $V_r\subseteq\mathbb{R}^p$ and an $n\ge \widetilde{n}_{r-1}$ be given
such that 
\begin{equation}
\xi_n\Bigl(f_{\theta_n}^{-1}\bigl(\overline{V}_r(\delta)\bigr)\Bigr)\,>\,\overline{\alpha}_r.
\label{eq2-add7}
\end{equation}
By property (\ref{eqC-s1-6}),  
$\delta^2/\kappa< \kappa\big/\bigl((c_r+1)^2\gamma^2\bigr)\le\kappa/\gamma^2\le1$, where the last inequality 
is obvious by the definitions of $\gamma$ and $\kappa$ in (\ref{eq2-2}) and (\ref{eq2-3}).
So $\delta<\sqrt{\kappa}$ and by Lemma \ref{lem2-4} and (\ref{eq2-add7})
\begin{equation}
f_{\theta_n}^\trp(x_{n+1})\,M^{-1}(\xi_n,\theta_n)\,f_{\theta_n}(x_{n+1})\,>\,
\Bigl(1-\bigl(1-\delta^2/\kappa\bigr)\overline{\alpha}_r\Bigr)^{-1}.
\label{eq2-add8}
\end{equation}
Next, we construct a particular basis $b_1,\ldots,b_r$ of the linear subspace $V_r$. 
From (\ref{eqC-13}) and (\ref{eq2-add7}) it follows that for all linear subspaces $V_t\subseteq\mathbb{R}^p$ of dimension $t\le r-1$ 
one has
\begin{equation}
\xi_n\Bigl(f_{\theta_n}^{-1}\bigl(\overline{V}_r(\delta)\setminus\overline{V}_t(\varepsilon_{r-1})\bigr)\Bigr)\,>\,
\overline{\alpha}_r-\alpha_{r-1}>0. \label{eqC-18}
\end{equation}
Note that by (\ref{eqC-s1-4}), in particular, the sets $R_1,\ldots,R_{K_r}$ cover ${\cal X}$. Thus (\ref{eqC-18})
implies that to any linear subspace $V\subseteq\mathbb{R}^p$ of dimension at most $r-1$ one can find some index $k\in\{1,\ldots,K_r\}$
such that $\xi_n\Bigl(f_{\theta_n}^{-1}\bigl(\overline{V}_r(\delta)\setminus\overline{V}(\varepsilon_{r-1})\bigr)\cap R_k)\Bigr)\,>\,
(\overline{\alpha}_r-\alpha_{r-1})/K_r$. Using this, one obtains inductively $r$ subsets $S_1,\ldots,S_r$
of ${\cal X}$ such that for all $j=1,\ldots,r$,
\begin{eqnarray}
&& S_j\,=\, f^{-1}_{\theta_n}\bigl(\overline{V}_r(\delta)\setminus\overline{W}_{j-1}(\varepsilon_{r-1})\bigr)\,\cap R_{k_j}
\ \mbox{ with some $k_j\in\{1,\ldots,K_r\}$,}\label{eqC-s2-1}\\
&& \xi_n(S_j)>(\overline{\alpha}_r-\alpha_{r-1})/K_r,\label{eqC-s2-2}
\end{eqnarray}
with particular linear subspaces $W_0,\ldots,W_{r-1}$ given by
\begin{equation} 
W_0\,=\,\{0\},\ \ W_t\,=\,{\rm span}\bigl\{\overline{f}_{\theta_n}(S_i,\xi_n)\,:\,1\le i\le t\bigr\}\ \mbox{for $1\le t\le r-1$,}
\label{eqC-s2-3}
\end{equation}
where $\overline{f}_{\theta_n}(S_i,\xi_n)$ denotes the average of $f_{\theta_n}(x)$ over $x\in S_i$ w.r.t. $\xi_n$ 
analogously to (\ref{eq2-add}). 
For each $j=1,\ldots,r$ by (\ref{eqC-s2-1}), firstly,   
$f_{\theta_n}(x)\in\overline{V}_r(\delta)$ for all $x\in S_j$ and hence also for the mean  
$\overline{f}_{\theta_n}(S_j,\xi_n)\in \overline{V}_r(\delta)$
since the set $\overline{V}_r(\delta)$ is convex. Secondly, $f_{\theta_n}(x)\not\in \overline{W}_{j-1}(\varepsilon_{r-1})$
for all $x\in S_j$, i.e.,
${\rm dist}(f_{\theta_n}(x),W_{j-1})>\varepsilon_{r-1}$ for all $x\in S_j$. Thirdly, $S_j\subseteq R_{k_j}$, hence
$\Vert f_{\theta_n}(x)-f_{\theta_n}(z)\Vert\le \varepsilon_{r-1}/2$ for all $x,z\in S_j$ which implies
 $\Vert f_{\theta_n}(x)-\overline{f}_{\theta_n}(S_j,\xi_n)\Vert\le\varepsilon_{r-1}/2$ for all $x\in S_j$. 
Using the inequality $\big|{\rm dist}(f_{\theta_n}(x),W_{j-1})-{\rm dist}\bigl(\overline{f}_{\theta_n}(S_j,\xi_n),W_{j-1}\bigr)\big|\le
\Vert f_{\theta_n}(x)-\overline{f}_{\theta_n}(S_j,\xi_n)\Vert$ one gets, choosing any $x\in S_j$, 
\[
{\rm dist}\bigl(\overline{f}_{\theta_n}(S_j,\xi_n),W_{j-1}\bigr)\,\ge\, {\rm dist}(f_{\theta_n}(x),W_{j-1})\,
-\,\Vert f_{\theta_n}(x)-\overline{f}_{\theta_n}(S_j,\xi_n)\Vert\,>\,\varepsilon_{r-1}-\varepsilon_{r-1}/2\,=\, 
\varepsilon_{r-1}/2.
\]
By (\ref{eqC-s2-3}) together with (M4) of the appendix  
the $p\times r$ matrix $F:=\bigl[\overline{f}_{\theta_n}(S_1,\xi_n),\ldots,\overline{f}_{\theta_n}(S_r,\xi_n)\bigr]$ satisfies 
\begin{equation}
\det\bigl(F^\trp F\bigr)>
\Bigl(\frac{\varepsilon_{r-1}^2}{4}\Bigr)^r. \label{eqC-20}
\end{equation}
For each $j=1,\ldots,r$, by $\overline{f}_{\theta_n}(S_j,\xi_n)\in\overline{V}_r(\delta)$, 
\begin{equation}
\overline{f}_{\theta_n}(S_j,\xi_n)\,=\,b_j + e_j\ \ \mbox{with }\ b_j\in V_r\ \mbox{and } e_j\in V_r^\perp,\ 
\Vert e_j\Vert\le\delta.\label{eqC-20a}
\end{equation}
Consider the $p\times r$ matrix $B:=\bigl[b_1,\ldots,b_r\bigr]$. Since 
$\Vert b_j\Vert\le\Vert\overline{f}_{\theta_n}(S_j,\xi_n)\Vert\le\gamma$ one has $F,B\in{\cal A}_r$.
By (\ref{eqC-20a}) $\Vert \overline{f}_{\theta_n}(S_j,\xi_n)-b_j\Vert\le\delta$, $1\le j\le r$, and hence, using property 
(\ref{eqC-s1-8}) of $\delta$ and (\ref{eqC-20}),
\[
\big|\det\bigl(F^\trp F\bigr)-\det\bigl(B^\trp B\bigr)\big|\,\le\,
\frac{1}{2}\Bigr(\frac{\varepsilon_{r-1}^2}{4}\Bigr)^r,\ \mbox{ and thus }\  
\det\bigl(B^\trp B\bigr)>\frac{1}{2}\Bigr(\frac{\varepsilon_{r-1}^2}{4}\Bigr)^r.
\]
In particular, $B\in{\cal A}_r^*$ 
and the vectors $b_1,\ldots,b_r$ are linearly independent and form thus a basis of the linear subspace 
$V_r$. Now suppose, contrary to the assertion of Step 2,  that $x_{n+1}\in f_{\theta_n}^{-1}\bigl(\overline{V}_r(\delta)\bigr)$.  
Then
\[
f_{\theta_n}(x_{n+1})\,=\,v+e\ \mbox{ for some }v\in V_r\ \mbox{and }e\in V_r^\perp, \ \Vert e\Vert\le\delta.
\] 
Since $b_1,\ldots,b_r$ constitute a basis of $V_r$ and $B=\bigl[b_1,\ldots,b_r\bigr]$,
one has \ $v\,=\,Bw$\ for some $w=(w_1,\ldots,w_r)^\trp\in\mathbb{R}^r$.  
In fact, $w$ is uniquely determined by $w=(B^\trp B)^{-1}B^\trp v$. Since 
$\Vert v\Vert\le\Vert f_{\theta_n}(x_{n+1})\Vert\le\gamma$ and $B\in{\cal A}_r^*$ one has, according to the definition
of $c_r$ in (\ref{eqC-s1-3}), that $\Vert w\Vert_1\le c_r$. Together with (\ref{eqC-20a}),
\begin{eqnarray*}
f_{\theta_n}(x_{n+1}) &=& \sum_{j=1}^rw_j b_j +e\,=\,\sum_{j=1}^rw_j\bigl(\overline{f}_{\theta_n}(S_j,\xi_n)-e_j\bigr) +e\\
& = & \sum_{j=1}^rw_j\overline{f}_{\theta_n}(S_j,\xi_n)-\sum_{j=1}^rw_je_j +e,\quad\mbox{and hence}
\end{eqnarray*}
\[
\Big\Vert f_{\theta_n}(x_{n+1})-\sum_{j=1}^rw_j\overline{f}_{\theta_n}(S_j,\xi_n)\Big\Vert\,\le\,\sum_{j=1}^r|w_j|\delta+\delta\,\le
\,(c_r+1)\delta.
\]
Define \ $a\,:=\,\sum_{j=1}^rw_j\overline{f}_{\theta_n}(S_j,\xi_n)$. Hence
\begin{equation}
\Vert f_{\theta_n}(x_{n+1})-a\Vert\,\le\,(c_r+1)\delta. \label{eqC-21a}
\end{equation}
Let $\eta:=(c_r+1)\,\delta\,\gamma/\kappa$. Then $\eta\in(\,0\,,\,1\,)$ by property 
(\ref{eqC-s1-6}) of $\delta$, 
and by (\ref{eqC-21a}) $\Vert f_{\theta_n}(x_{n+1})-a\Vert\le\eta\kappa/\gamma$. So, by Lemma \ref{lem2-2},
\begin{equation}
f_{\theta_n}^\trp(x_{n+1})\,M^{-1}(\xi_n,\theta_n)\,f_{\theta_n}(x_{n+1})\,\le\,
\frac{1}{(1-\eta)^2}\,a^\trp M^{-1}(\xi_n,\theta_n)\,a.\label{eqC-22}
\end{equation}
Observing that $b\mapsto \bigl(b^\trp M^{-1}(\xi_n,\theta_n)\,b\bigr)^{1/2}$, $b\in\mathbb{R}^p$, is a norm on 
$\mathbb{R}^p$ and using the definition of the vector $a$, 
\begin{equation}
\bigl(a^\trp M^{-1}(\xi_n,\theta_n)\,a\bigr)^{1/2}\,\le\,
\sum_{j=1}^r|w_j|\,\Bigl(\overline{f}_{\theta_n}^\trp(S_j,\xi_n)\,M^{-1}(\xi_n,\theta_n)\,\overline{f}_{\theta_n}(S_j,\xi_n)\Bigr)^{1/2}.\label{eqC-22a}
\end{equation}
For each $j=1,\ldots,r$, one gets by (M1) and (M2) of the appendix, 
where the sums below are taken over $x\in S_j\cap\,{\rm supp}(\xi_n)$,
\begin{eqnarray*}
&&\overline{f}_{\theta_n}^\trp(S_j,\xi_n)\,M^{-1}(\xi_n,\theta_n)\,\overline{f}_{\theta_n}(S_j,\xi_n)\\
&&\le\ \Bigl(\sum_x\frac{\xi_n(x)}{\xi_n(S_j)}f_{\theta_n}(x)\Bigr)^\trp\Bigl(\xi_n(S_j)\sum_x
\frac{\xi_n(x)}{\xi_n(S_j)}f_{\theta_n}(x)f_{\theta_n}^\trp(x)\Bigr)^-\Bigl(\sum_x\frac{\xi_n(x)}{\xi_n(S_j)}f_{\theta_n}(x)\Bigr)\\
&&\le\ 1\big/\xi_n(S_j)\,<\,\frac{K_r}{(\overline{\alpha}_r-\alpha_{r-1})},
\end{eqnarray*}
where the last inequality is due to (\ref{eqC-s2-2}). Hence by (\ref{eqC-22a}) and by $\Vert w\Vert_1\le c_r$, 
\[
\bigl(a^\trp M^{-1}(\xi_n,\theta_n)\,a\bigr)^{1/2}\,\le\,
\Bigl(\frac{K_r}{\overline{\alpha}_r-\alpha_{r-1}}\Bigr)^{1/2}\,c_r,
\]
and together with (\ref{eqC-22}) one gets
\begin{eqnarray}
&&f_{\theta_n}^\trp(x_{n+1})\,M^{-1}(\xi_n,\theta_n)\,f_{\theta_n}(x_{n+1})\,\le\, 
\frac{K_rc_r^2}{(1-\eta)^2(\overline{\alpha}_r-\alpha_{r-1})},
\label{eqC-s2-5}\\
&&\mbox{where, as before,}\ \eta\,=\,(c_r+1)\delta\gamma/\kappa.\nonumber
\end{eqnarray}
Observing that the r.h.s. of (\ref{eqC-s2-5}) equals the reciprocal of the l.h.s. of 
(\ref{eqC-s1-7}), it follows from (\ref{eqC-s1-7}) that
\[
\frac{K_rc_r^2}{(1-\eta)^2(\overline{\alpha}_r-\alpha_{r-1})}\,<\,
\Bigl(1-\bigl(1-\delta^2/\kappa\bigr)\,\overline{\alpha}_r\Bigr)^{-1},
\]
and hence by (\ref{eqC-s2-5})
\[
f_{\theta_n}^\trp(x_{n+1})\,M^{-1}(\xi_n,\theta_n)\,f_{\theta_n}(x_{n+1})\,<\,\Bigl(1-\bigl(1-\delta^2/\kappa\bigr)\,\overline{\alpha}_r\Bigr)^{-1},
\]
which is a contradiction to (\ref{eq2-add8}) derived above. So our supposition that $x_{n+1}\in f_{\theta_n}^{-1}\bigl(\overline{V}_r(\delta)\bigr)$
was wrong. Hence the result of Step 2 follows.\\[.5ex] 
\underline{Step 3.} \ 
For $\delta$ from the previous Steps 1 and 2, let $m_0(\delta)$ and $\delta'>0$ be chosen according to (T),
where we may assume that $m_0(\delta)\ge \widetilde{n}_{r-1}$. 
Recall that  $0<\overline{\alpha}_r<1$ according to (\ref{eqC-s1-5}). 
Choose an  $\alpha_r$ such that $\overline{\alpha}_r<\alpha_r<1$.
Let $V_r$ be any $r$-dimensional linear subspace of $\mathbb{R}^p$. Consider the set 
\[
 S\,=\,\bigcap_{\ell\ge m_0(\delta)}f_{\theta_\ell}^{-1}\bigl(\overline{V}_r(\delta)\bigr).
\]
By the result of Step 2 and by \ $S\subseteq f_{\theta_n}^{-1}\bigl(\overline{V}_r(\delta)\bigr)$ for all $n\ge m_0(\delta)$, we have:
\[
\mbox{If $n\ge m_0(\delta)$ and $\xi_n(S)>\overline{\alpha}_r$ \ then \ $x_{n+1}\not\in S$.}
\]
So the sequence 
\[
\beta_n:=\xi_n(S), \ n\ge m_0(\delta),\  
\mbox{ along with $\beta:=\overline{\alpha}_r$ and $\widetilde{\beta}:=\alpha_r$}
\] 
satisfies the assumptions of Lemma \ref{lem2-1}, and hence by that lemma
$\beta_n\le \alpha_r$ for all $n\ge m_1=m_1\bigl(\beta,\widetilde{\beta},m_0(\delta)\bigr)$. 
Since $m_1$ does not depend on the particular choice of $V_r$ we have thus obtained that  
for all linear subspaces $V_r\subseteq\mathbb{R}^p$ of dimension $r$ one has
\[
 \xi_n\Bigl(\bigcap_{\ell\ge m_0(\delta)}f_{\theta_\ell}^{-1}\bigl(\overline{V}_r(\delta)\bigr)\Bigr)\,\le\,\alpha_r\quad
\mbox{for all $n\ge m_1$.}
\]
According to (T) we have $\delta'>0$ such that  
$f_{\theta_n}^{-1}\bigl(\overline{V}_r(\delta')\bigr)\subseteq \bigcap_{\ell\ge m_0(\delta)}f_{\theta_\ell}^{-1}\bigl(\overline{V}_r(\delta)\bigr)$
for all $n\ge m_1$ and all linear subspaces $V_r$ of dimension $r$. Hence 
$\xi_n\Bigl(f_{\theta_n}^{-1}\bigl(\overline{V}_r(\delta')\bigr)\Bigr)\le\alpha_r$ for all $n\ge m_1$ and all linear subspaces $V_r$ of dimension $r$,
which is statement ${\rm S}(r)$ with $\widetilde{n}_r:=m_1$, $\varepsilon_r:=\delta'$, and $\alpha_r$ as obtained.
So the induction step has been completed. 
\eop 

\begin{corollary}\quad\label{cor2-1}
(i) If (T) is satisfied  then the asymptotic nonsingularity (\ref{eq2-1a}) holds.\\
(ii) If (GLM) is satisfied then the stronger asymptotic nonsingularity (\ref{eq2-1b}) holds.
\end{corollary}
   
\noindent{\bf Proof.} \  
Using a well-known representation of the smallest eigenvalue of a symmetric matrix we can write
\begin{equation}
\lambda_{\rm\scriptsize min}\bigl(M(\xi_n,\theta)\bigr)\,=\,\min_{\Vert c\Vert=1}c^\trp M(\xi_n,\theta)\,c\,=\,
\min_{\Vert c\Vert=1}\sum_{x\in{\rm supp}(\xi_n)}\xi_n(x)\,\bigl(c^\trp f_\theta(x)\bigr)^2. \label{eq2-5}
\end{equation}
For any $c\in\mathbb{R}^p$, $\Vert c\Vert=1$, we denote by $V_{p-1,c}$ the $(p-1)$-dimensional subspace of $\mathbb{R}^p$ given by
$V_{p-1,c}\,=\,\bigl\{a\in\mathbb{R}^p\,:\,c^\trp a=0\bigr\}$. 
Assume (T).
Let $n_0\ge n_{\rm\scriptsize st}$, $\varepsilon>0$, and $\alpha\in(\,0\,,\,1\,)$
be chosen according to Theorem \ref{theo2-1}. Then by the theorem, observing that \ 
$f^{-1}_{\theta_n}\bigl(\overline{V}_{p-1,c}(\varepsilon)\bigr)\,=\,\bigl\{x\in{\cal X}\,:\,|c^\trp f_{\theta_n}(x)|\le\varepsilon\bigr\}$,
we have for all $n\ge n_0$ and all $c$ with $\Vert c\Vert =1$,
\[ 
\xi_n\Bigl(\bigl\{x\in{\cal X}\,:\,|c^\trp f_{\theta_n}(x)|\le\varepsilon\bigr\}\Bigr)\,\le\alpha.
\] 
Denote 
$S_{n,c}\,=\,\bigl\{x\in{\cal X}\,:\,|c^\trp f_{\theta_n}(x)|>\varepsilon\bigr\}$.
Then for all $n\ge n_0$ and all $c\in\mathbb{R}^p$ with $\Vert c\Vert =1$ one has $\xi_n(S_{n,c})\ge 1-\alpha$  and 
hence
\[
\sum_{x\in{\rm supp}(\xi_n)}\xi_n(x)\,\bigl(c^\trp f_{\theta_n}(x)\bigr)^2\,\ge\,
\sum_{x\in S_{n,c}\cap\,{\rm supp}(\xi_n)}\xi_n(x)\,\bigl(c^\trp f_{\theta_n}(x)\bigr)^2\,\ge\,
\varepsilon^2(1-\alpha)\,=:\widetilde{\lambda}_0>0,
\]
and together with (\ref{eq2-5}), $\lambda_{\rm\scriptsize min}\bigl(M(\xi_n,\theta_n)\bigr)\ge \widetilde{\lambda}_0$
for all $n\ge n_0$.  
So, in the case  $n_0>n_{\rm\scriptsize st}$, a positive real constant is given by 
$\lambda_0:=\min\{\widetilde{\lambda}_0,\lambda_{\rm\scriptsize min}\bigl(M(\xi_n,\theta_n)\bigr):\,n_{\rm\scriptsize st}\le n<n_0\}$.
In the case $n_0=n_{\rm\scriptsize st}$ choose   
$\lambda_0:=\widetilde{\lambda}_0$. In any case,  with that constant $\lambda_0>0$ the asymptotic nonsingularity
(\ref{eq2-1a}) holds. Now assume (GLM). By Lemma \ref{lem2-4} (T) is satisfied as well and hence, as already
proved, the asymptotic nonsingularity (\ref{eq2-1a}) holds with some $\lambda_0>0$. 
For all $x\in{\cal X}$ and all $\theta\in\Theta$ one has $f_\theta(x)\,f_\theta^\trp(x)=
\psi^2(x,\theta)\,f(x)\,f^\trp(x)$. Consider the positive real numbers $\psi_{\rm\scriptsize min}$ and
$\psi_{\rm\scriptsize max}$ from (\ref{eq2-add10}). Then, for all $x\in{\cal X}$ and $\theta,\theta'\in\Theta$ trivially
$(\psi_{\rm\scriptsize min}/\psi_{\rm\scriptsize max})^2\psi^2(x,\theta')\le \psi^2(x,\theta)$. Hence for any design $\xi$
one has $(\psi_{\rm\scriptsize min}/\psi_{\rm\scriptsize max})^2M(\xi,\theta')\le M(\xi,\theta)$ for all $\theta,\theta'\in\Theta$.
In particular, one has
$(\psi_{\rm\scriptsize min}/\psi_{\rm\scriptsize max})^2M(\xi_n,\theta_n)\le M(\xi_n,\theta)$ for all $n\ge n_{\rm\scriptsize st}$
and $\theta\in\Theta$. It follows that
\[
\lambda_{\rm\scriptsize min}\bigl(M(\xi_n,\theta)\bigr)\,\ge\,(\psi_{\rm\scriptsize min}/\psi_{\rm\scriptsize max})^2\,
\lambda_{\rm\scriptsize min}\bigl(M(\xi_n,\theta_n)\bigr) \ge 
(\psi_{\rm\scriptsize min}/\psi_{\rm\scriptsize max})^2\lambda_0\quad\mbox{for all $n\ge n_{\rm\scriptsize st}$, $\theta\in\Theta$.} 
\]
So the stronger asymptotic nonsingularity (\ref{eq2-1b}) holds with $(\psi_{\rm\scriptsize min}/\psi_{\rm\scriptsize max})^2\lambda_0$
instead of $\lambda_0$.
\eop    

As a further consequence from Theorem \ref{theo2-1} and Corollory \ref{cor2-1} we can derive a convergence result  
as in Pronzato \cite{Pronzato}, Lemma 2 and Theorem 2, 
and Freise \cite{Fritjof}, Lemma 18. If the sequence
of parameter points $\theta_n$ converges to some parameter point $\overline{\theta}\in\Theta$ 
then the design sequence $\xi_n$ is asymptotically locally D-optimal at $\overline{\theta}$,
in the sense that the sequence of information matrices $M(\xi_n,\theta_n)$ converges to the 
information matrix $M(\xi^*_{\overline{\theta}},\overline{\theta})$ of a locally D-optimal design $\xi^*_{\overline{\theta}}$
at $\overline{\theta}$. 
For later reference  (see Section 3), the next theorem states the convergence of the information matrices
$M(\xi_n,\theta_n')$ to $M(\xi^*_{\overline{\theta}},\overline{\theta})$ for {\em any} sequence $\theta_n'\in\Theta$ 
converging to $\overline{\theta}$, provided that the sequence $\theta_n$ employed by the algorithm converges to 
$\overline{\theta}$. 
Of course, in the linear model case, $f_\theta=f$ identical for all $\theta\in\Theta$, we retrieve 
the classical result of Wynn \cite{Wynn}, Theorem 1. 

\begin{theorem}\quad\label{theo2-2}\\
If \ $\lim_{n\to\infty}\theta_n=\overline{\theta}$ \ for some $\overline{\theta}\in\Theta$ then for every sequence $\theta_n'\in\Theta$,
$n\ge n_{\rm\scriptsize st}$, such that\\ 
$\lim_{n\to\infty}\theta_n'=\overline{\theta}$ one has 
\[
\lim_{n\to\infty}M(\xi_n,\theta_n')\,=\,M(\xi^*_{\overline{\theta}},\overline{\theta}),
\] 
where $\xi^*_{\overline{\theta}}$ denotes a locally D-optimal design at $\overline{\theta}$, i.e., $\xi^*_{\overline{\theta}}$ maximizes 
$\det\bigl(M(\xi,\overline{\theta})\bigr)$ over the set of all designs $\xi$. 
\end{theorem}

\noindent
{\bf Proof.} \ 
The matrix-valued function \ $(x,\theta)\longmapsto f_\theta(x)\,f_\theta^\trp(x)$ \ is uniformly continuous on its compact domain 
${\cal X}\times\Theta$. So, for any sequence $\theta_n'\in\Theta$ converging to $\overline{\theta}$, 
observing that $M(\xi_n,\theta)=\frac{1}{n}\sum_{i=1}^nf_\theta(x_i)\,f_\theta^\trp(x_i)$ for all $n\ge n_{\rm\scriptsize st}$
and $\theta\in\Theta$,
\begin{eqnarray*}
&&\big\Vert M(\xi_n,\theta_n')-M(\xi_n,\overline{\theta})\big\Vert\,\le\,
\frac{1}{n}\sum_{i=1}^n
\big\Vert f_{\theta_n'}(x_i)\,f_{\theta_n'}^\trp(x_i)-  f_{\overline{\theta}}(x_i)\,f_{\overline{\theta}}^\trp(x_i) \big\Vert\\
&&\,\le\,\max_{x\in{\cal X}}\big\Vert f_{\theta_n'}(x)\,f_{\theta_n'}^\trp(x)-  f_{\overline{\theta}}(x)\,f_{\overline{\theta}}^\trp(x) \big\Vert\ 
\longrightarrow 0\ \ \mbox{ as $n\to\infty$}.  
\end{eqnarray*}    
Hence 
\begin{equation}
\big\Vert M(\xi_n,\theta_n')-M(\xi_n,\overline{\theta})\big\Vert\ \longrightarrow 0\ \ \mbox{ as $n\to\infty$}, \label{eq2-6a}
\end{equation}
and, in particular, 
\begin{equation}
\big\Vert M(\xi_n,\theta_n)-M(\xi_n,\overline{\theta})\big\Vert\ \longrightarrow 0 \ \mbox{ as $n\to\infty$.}
\label{eq2-6c}
\end{equation} 
Consider $\gamma$ from (\ref{eq2-2}). For any design $\xi$ and any $\theta\in\Theta$ we have
\[
{\rm tr}\bigl(M(\xi,\theta)\bigr)=\sum_{x\in{\rm\scriptsize supp}(\xi)}\xi(x)\,\Vert\,f_\theta(x)\Vert^2\le\gamma^2. 
\] 
By Lemma \ref{lem2-4} and Corollary \ref{cor2-1} there is a $\lambda_0>0$ satisfying (\ref{eq2-1a}). 
Let ${\cal A}$ be the set of all nonnegative definite $p\times p$ martices $A$ such that $\lambda_{\rm\scriptsize min}(A)\ge\lambda_0/2$
and ${\rm tr}(A)\le\gamma^2$. Clearly, ${\cal A}$ is compact and $M(\xi_n,\theta_n)\in{\cal A}$ for all $n\ge n_{\rm\scriptsize st}$, 
and by (\ref{eq2-6c}) there is an $m_0\ge n_{\rm\scriptsize st}$ such that $M(\xi_n,\overline{\theta})\in{\cal A}$ for all $n\ge m_0$.  
Define a real-valued function $H$ on ${\cal X}\times \Theta\times {\cal A}$ by
\[
H(x,\theta,A)\,=\,f^\trp_\theta(x)\,A^{-1} f_\theta(x),
\]
which is continuous and hence uniformly continuous on its compact domain ${\cal X}\times \Theta\times {\cal A}$.
So, together with (\ref{eq2-6c}),
\begin{equation}
\max_{x\in{\cal X}}\big|H\bigl(x,\theta_n,M(\xi_n,\theta_n)\bigr) - 
H\bigl(x,\overline{\theta},M(\xi_n,\overline{\theta})\bigr)\big|\,\longrightarrow 0\ \mbox{ as $n\to\infty$.}\label{eq2-6b}
\end{equation}
In what follows let an $\varepsilon\in(\,0\,,\,1\,)$ be given. By (\ref{eq2-6b}) and by the definition of the function $H$,
there is an $n_0\ge m_0$ such that
\[
\big|f_{\theta_n}^\trp(x)M^{-1}(\xi_n,\theta_n)\,f_{\theta_n}(x) - 
f_{\overline{\theta}}^\trp(x)M^{-1}(\xi_n,\overline{\theta})\,f_{\overline{\theta}}(x)\big|\ \le\varepsilon/2
\ \mbox{ for all $x\in{\cal X}$ and all $n\ge n_0$.}
\] 
This yields, in particular, \\[.5ex]
{\bf(i)} $\forall n\ge n_0$\,:\ \ 
$\displaystyle f^\trp_{\overline{\theta}}(x_{n+1})\,M^{-1}(\xi_n,\overline{\theta})\,f_{\overline{\theta}}(x_{n+1})\,\ge\,
\max_{x\in{\cal X}}f^\trp_{\overline{\theta}}(x)\,M^{-1}(\xi_n,\overline{\theta})\,
f_{\overline{\theta}}(x)\,-\,\varepsilon\,\ge\,p-\varepsilon$,\\[.5ex]
since for all $n\ge n_0$, denoting \ $x_n^*=\arg\max_{x\in{\cal X}}f_{\overline{\theta}}^\trp(x)\,M^{-1}(\xi_n,\overline{\theta})\,f_{\overline{\theta}}(x)$,
\begin{eqnarray*}
&&f_{\overline{\theta}}^\trp(x_{n+1})\,M^{-1}(\xi_n,\overline{\theta})\,f_{\overline{\theta}}(x_{n+1})\,\ge 
f_{\theta_n}^\trp(x_{n+1})\,M^{-1}(\xi_n,\theta_n)\,f_{\theta_n}(x_{n+1})\,-\,\varepsilon/2\\
&&\ge\,f_{\theta_n}^\trp(x^*_n)\,M^{-1}(\xi_n,\theta_n)\,f_{\theta_n}(x^*_n)\,-\,\varepsilon/2  
\,\ge\,f_{\overline{\theta}}^\trp(x^*_n)\,M^{-1}(\xi_n,\overline{\theta})\,f_{\overline{\theta}}(x^*_n)\,-\varepsilon.
\end{eqnarray*}
The second inequality in (i) is well-known from the Kiefer-Wolfowitz Equivalence Theorem. The rest of the proof 
employs the arguments of Pronzato \cite{Pronzato} in the proof of Lemma 3 of that paper. For convenience we report here
the main steps labelled below by (ii) - (v).\\[.5ex]
{\bf(ii)} One can choose $n_1\ge n_0$ such that for all $n\ge n_1$
\[ 
\log\det\bigl(M(\xi_{n+1},\overline{\theta})\bigr)\,-\,
\log\det\bigl(M(\xi_n,\overline{\theta})\bigr)\,\ge\,-\varepsilon.
\]
To see this we note that  \ $(n+1)M(\xi_{n+1},\overline{\theta})=nM(\xi_n,\overline{\theta}) + f_{\overline{\theta}}(x_{n+1})
\,f^\trp_{\overline{\theta}}(x_{n+1})$ and by a well-known formula of determinants,
\begin{eqnarray}
&&\log\det\bigl(M(\xi_{n+1},\overline{\theta})\bigr)\,-\,\log\det\bigl(M(\xi_n,\overline{\theta})\bigr)\,= \nonumber\\
&&\phantom{xxx}\log\Bigl(1+\frac{1}{n}f^\trp_{\overline{\theta}}(x_{n+1})\,M^{-1}(\xi_n,\overline{\theta})\,f_{\overline{\theta}}(x_{n+1})\Bigr)
\,-\,p\log\Bigl(1+\frac{1}{n}\Bigr). \label{eq2-6}
\end{eqnarray}
By (i) for $n\ge n_0$ the expression (\ref{eq2-6}) is greater than or equal to  
\[
\log\Bigl(1+\frac{p-\varepsilon}{n}\Bigr)-p\log\Bigl(1+\frac{1}{n}\Bigr)\,=\,\log\Bigl(\frac{1+(p-\varepsilon)/n}
{\bigl(1+1/n\bigr)^p}\Bigr)\,=\log\Bigl(\frac{1+(p-\varepsilon)/n}{1+(p+c_n)/n}\Bigr),
\]
where we have used that $(1+1/n)^p=1+(p+c_n)/n$ with $c_n>0$, $c_n\to0$ as $n\to\infty$.
Choose $n_1\ge n_0$ such that $c_n\le(p-\varepsilon)\varepsilon$ for all $n\ge n_1$. Then for all $n\ge n_1$,
\begin{eqnarray*}
&&\log\Bigl(\frac{1+(p-\varepsilon)/n}{1+(p+c_n)/n}\Bigr)\,\ge\,
\log\Bigl(\frac{1+(p-\varepsilon)/n}{1+\bigl(p+(p-\varepsilon)\varepsilon)/n}\Bigr)\,\\
&&\,\ge\,-\frac{1}{1+(p-\varepsilon)/n}\,\frac{p+(p-\varepsilon)\varepsilon-(p-\varepsilon)}{n}
\,=\,-\frac{\varepsilon(1+p-\varepsilon)}{n+p-\varepsilon}\,\ge\,-\varepsilon.
\end{eqnarray*}
{\bf(iii)} One can choose $n_2\ge n_1$ such that for all $n\ge n_2$ 
\[
\log\Bigl(1+\frac{p+\varepsilon}{n}\Bigr)\,-\,p\log\Bigl(1+\frac{1}{n}\Bigr)\,\ge\frac{\varepsilon}{2n}.
\]
This follows from $(1+1/n)^p=1+(p+c_n)/n$  and by choosing $n_2\ge n_1$
such that $c_n\le\varepsilon\bigl(1-\frac{n+p+\varepsilon}{2n}\bigr)$ for all $n\ge n_2$. \\[.5ex]
{\bf(iv)} Denote 
$\Psi^*:=\log\det\bigl(M(\xi^*_{\overline{\theta}},\overline{\theta})\bigr)$. 
If $n\ge n_2$ and \ $\log\det\bigl(M(\xi_n,\overline{\theta})\bigr)\,\le\ \Psi^*-2\varepsilon$ \ then
\[
\log\det\bigl(M(\xi_{n+1},\overline{\theta})\bigr)\,-\,
\log\det\bigl(M(\xi_n,\overline{\theta})\bigr)\,\ge\,\frac{\varepsilon}{2n}.
\]
This can be seen as follows. By the gradient inequality for the concave criterion $\log\det(\,\cdot\,)$, 
\begin{eqnarray*}
&&2\varepsilon\,\le\,\Psi^*\,-\,\log\det\bigl(M(\xi_n,\overline{\theta})\bigr)\,\le
\,\max_{x\in{\cal X}}f_{\overline{\theta}}^\trp(x)M^{-1}(\xi_n,\overline{\theta})\,
f_{\overline{\theta}}(x)\,-p\\
&&\,\le\,f_{\overline{\theta}}^\trp(x_{n+1})M^{-1}(\xi_n,\overline{\theta})\,
f_{\overline{\theta}}(x_{n+1})\,+\varepsilon -p,
\end{eqnarray*}
where the last inequality comes from (i). Hence it follows that 
$f_{\overline{\theta}}^\trp(x_{n+1})M^{-1}(\xi_n,\overline{\theta})\,
f_{\overline{\theta}}(x_{n+1})\ge p+\varepsilon$ and together with (\ref{eq2-6}) one gets
\[
\log\det\bigl(M(\xi_{n+1},\overline{\theta})\bigr)\,-\,\log\det\bigl(M(\xi_n,\overline{\theta})\bigr)\,\ge\,
\log\Bigl(1+\frac{p+\varepsilon}{n}\Bigr)\,-\,p\log\Bigl(1+\frac{1}{n}\Bigr)\,\ge\,\frac{\varepsilon}{2n},
\] 
where the last inequality comes from (iii).\\[.5ex]
{\bf(v)} One can choose $n_3\ge n_2$ such that for all $n\ge n_3$
\[
\log\det\bigl(M(\xi_n,\overline{\theta})\bigr)\,>\,\Psi^*-3\varepsilon.
\]
To see this, note that by (iv) there is some $n_3\ge n_2$ such that
\ $\log\det\bigl(M(\xi_{n_3},\overline{\theta})\bigr)>\Psi^*-2\varepsilon$, since otherwise (iv) would yield that
$\log\det\bigl(M(\xi_n,\overline{\theta})\bigr)\longrightarrow\infty$ as $n\to\infty$, which is a contradiction.
By (ii) and (iv), the sequence \ $a_n:=\log\det\bigl(M(\xi_n,\overline{\theta})\bigr)$, $n\ge n_3$,  has the following properties. 
\[
a_{n_3}>\Psi^*-2\varepsilon\,;\ \ a_{n+1}-a_n\ge -\varepsilon\ \ \forall\ n\ge n_3;\ \ 
a_{n+1}-a_n>0\ \mbox{ if $a_n\le\Psi^*-2\varepsilon$}.
\] 
Thus, obviously, $a_n>\Psi^*-3\varepsilon$ for all $n\ge n_3$, which is (v).

From (v) we get
\[
\liminf_{n\to\infty} \log\det\bigl(M(\xi_n,\overline{\theta})\bigr)\,\ge\,\Psi^*-3\varepsilon.
\]
Since $\varepsilon\in(\,0\,,\,1\,)$ was arbitrary we get  
\ $\liminf_{n\to\infty} \log\det\bigl(M(\xi_n,\overline{\theta})\bigr)\,\ge\,\Psi^*$ \ and hence
\ $\lim_{n\to\infty} \log\det\bigl(M(\xi_n,\overline{\theta})\bigr)\,=\,\Psi^*$.
This implies  \ $\lim_{n\to\infty}M(\xi_n,\overline{\theta})\,=\,M(\xi^*_{\overline{\theta}},\overline{\theta})$,
since by strict concavity of the criterion $\log\det(\,\cdot\,)$ the information matrix at $\overline{\theta}$  
of a locally D-optimal design at $\overline{\theta}$ is unique. 
That is, denoting by $\Xi$ the set of all designs and  
${\cal M}_{\overline{\theta}}:=\bigl\{M(\xi,\overline{\theta})\,:\,\mbox{$\xi\in\Xi$ }\bigr\}$,
the set of all information matrices of designs at $\overline{\theta}$, the information matrix
$M^*=M(\xi^*_{\overline{\theta}},\overline{\theta})$ is the unique point in  ${\cal M}_{\overline{\theta}}$ 
such that \ $\log\det(M^*)\,=\max_{M\in{\cal M}_{\overline{\theta}}}\log\det(M)$. So for any $\delta>0$ one has by
compactness and continuity
\[
\sup\bigl\{\log\det(M)\,:\,M\in{\cal M}_{\overline{\theta}},\ \Vert M-M^*\Vert\ge\delta\bigr\}\ <\,\log\det(M^*). 
\]
So, $\lim_{n\to\infty}\log\det\bigl(M(\xi_n,\overline{\theta}))\,=\log\det(M^*)$ 
implies $\lim_{n\to\infty}M(\xi_n,\overline{\theta})\,=M^*$.
If $\theta_n'\in\Theta$ is any sequence converging to $\overline{\theta}$ then by (\ref{eq2-6a})
$\lim_{n\to\infty}M(\xi_n,\theta_n')\,=M^*$.
\eop

\section{Adaptive Wynn-algorithm in univariate GLM}
\setcounter{equation}{0}
Now we focus on the adaptive character of the algorithm. The sequence of parameter points $\theta_n$, $n\ge n_{\rm\scriptsize st}$,
employed is given by parameter estimates based on the data available at the current stage $n$, which are the design points $x_1,\ldots,x_n$
{\em and} the observed values $y_1,\ldots,y_n$ of a univariate response variable.
We assume a (nonlinear) regression model with expected univariate responses $\mu(x,\theta)$, where $x\in{\cal X}$ and $\theta\in\Theta$.
The function $\mu\,:\,{\cal X}\times\Theta\longrightarrow\mathbb{R}$ is assumed to be continuous and, 
as in the previous sections, the experimental region
${\cal X}$ and the parameter space $\Theta$ are compact metric spaces.
Again, for the algorithm we assume a family $f_\theta$, $\theta\in\Theta$, 
of $\mathbb{R}^p$-valued functions on ${\cal X}$ defining the information matrices of designs by (\ref{eq1-0})
and having the properties that for each $\theta\in\Theta$ the image $f_\theta({\cal X})$ spans $\mathbb{R}^p$, and the function
$(x,\theta)\longmapsto f_\theta(x)$ is continuous on ${\cal X}\times\Theta$.  The adaptive Wynn-algorithm 
sequentially generates data $x_1,y_1,\ldots,x_n,y_n,\ldots$ where $y_i$ is the observed (univariate) response
at the design point $x_i$ ($i=1,2,\ldots,n,\ldots$) and 
the employed sequence $\theta_n$, $n\ge n_{\rm\scriptsize st}$, is given by adaptive 
parameter estimates, \ $\widehat{\theta}_n=\widehat{\theta}_n(x_1,y_1,\ldots,x_n,y_n)$, $n\ge n_{\rm\scriptsize st}$. 
In particular, the values $y_i$ of the response variable as well as the generated values $x_i$
of the design variable are random and hence they are modelled by random variables $Y_i$ and $X_i$.  
The sequential and adaptive character of the data is caught by the `adaptive regression model' formulated and discussed in 
Subsection 3.1 below. 
For theoretical investigations on consistency or asymptotic distribution of estimators it 
will be convenient to distinguish between the true (but unknown) parameter point $\overline{\theta}$ 
and any possible parameter point $\theta\in\Theta$ to be considered.    
So throughout this section, $\overline{\theta}\in\Theta$ denotes the fixed true parameter point governing the 
random variables. 

\subsection{Adaptive regression model.}
An appropriate model for the adaptive character of the  
sequences of random variables $X_i$ and $Y_i$, $i\in\mathbb{N}$, is provided by the  following assumptions 
(A1) and (A2), 
cp. Lai \cite{Lai}, Sec.~1,  or Chen, Hu, and Ying \cite{Chen-Hu-Ying}, Sec.~3. Note that all the random variables 
are defined on some probability space $(\Omega,{\cal F},\mathbb{P}_{\overline{\theta}})$, 
where $\Omega$ is a nonempty set, ${\cal F}$
is a sigma-field of subsets of $\Omega$, and $\mathbb{P}_{\overline{\theta}}$ is a probability measure on ${\cal F}$ 
corresponding to the true parameter point $\overline{\theta}$. 
\begin{itemize}
\item[{\bf(A1)}]
There is given a nondecreasing sequence of sub-sigma-fields of ${\cal F}$, 
\ ${\cal F}_0\subseteq{\cal F}_1\subseteq\,\ldots\,\subseteq{\cal F}_n\subseteq,\ldots$ 
such that  for each $i\in\mathbb{N}$ the random variable $X_i$ is ${\cal F}_{i-1}$-measurable
and the random variable $Y_i$ is ${\cal F}_i$-measurable.
\item[{\bf(A2)}]
$Y_i\,=\,\mu(X_i,\overline{\theta})\,+\,e_i$ \ with real-valued 
square integrable random errors $e_i$ such that\\ 
${\rm E}\bigl(e_i\,\big|\,{\cal F}_{i-1}\bigr)\,=0 \ \mbox{\,a.s.}$\  for all $i\in\mathbb{N}$, and \ 
$\sup_{i\in\mathbb{N}}{\rm E}\bigl(e_i^2\,\big|\,{\cal F}_{i-1}\bigr)\,<\infty \ \ \mbox{a.s.}$
\end{itemize}
As an illustration of the sub-sigma-fields ${\cal F}_i$, $i\in\mathbb{N}_0$, suppose that 
the starting design $\xi_{n_{\rm\tiny st}}$ of the algorithm was chosen deterministically, i.e., 
$X_1,\ldots,X_{n_{\rm\tiny st}}$ are constants, and suppose further that for all $n\ge n_{\rm\scriptsize st}$ there is no ambiguity in chosing
the maximizer $x_{n+1}$ in (\ref{eq1-2}) given the values of $X_1,Y_1,\ldots,X_n,Y_n$ and thus given the value of $\widehat{\theta}_n$. 
Then for all $n\ge n_{\rm\scriptsize st}$
the random variable $X_{n+1}$ is a function of $Y_1,\ldots,Y_n$. So one can employ the particular sigma-fields 
\ ${\cal F}_n=\sigma(Y_1,\ldots,Y_n)$ generated by the random variables $Y_1,\ldots,Y_n$, for all $n\ge1$,
and ${\cal F}_0$ the minimal sigma-field in $\Omega$.    
We note that no further relation is assumed so far between the mean response function $\mu$ and the family of functions
$f_\theta$, $\theta\in\Theta$, of the algorithm, whereas a particular relation will be employed in the next subsection.

The following lemma presents some auxiliary asymptotic results derived from martingale limit theorems.
If $W_n$, $n\in \mathbb{N}$, is a sequence of $\mathbb{R}^k$-valued random variables and $W$ is an $\mathbb{R}^k$-valued 
random variable, the notation \ $W_n\asto\,W$ \ stands for almost sure convergence
of the sequence $W_n$ to $W$ (as $n\to\infty$). For real-valued $W_n$ we will also use the notation  
$W_n\asto\,\infty$ \  for indicating almost sure convergence (or `divergence') to infinity. 

\begin{lemma}\quad\label{lem3-1}
Under(A1) and (A2) the following (a), (b), and (c) hold.\\ 
{\bf(a)} \ 
$\displaystyle\limsup_{n\to\infty}\,\frac{1}{n}\,\sum_{i=1}^n|e_i|\ <\,\infty\ \ \mbox{a.s.}$
\\[1ex]
{\bf(b)} \ Let $Z_i$, $i\in\mathbb{N}$, be a sequence of real-valued square integrable random variables such that  
$Z_i$ is ${\cal F}_{i-1}$-measurable for all $i\in\mathbb{N}$ and \ $\sup_{i\in\mathbb{N}}|Z_i|\,<\infty\,\mbox{ a.s.}$ 
Then 
\[
\frac{1}{n}\sum_{i=1}^nZ_ie_i\,\asto\,0\,.
\]
{\bf(c)} \ Let \ $h\,:\,{\cal X}\times\Theta\longrightarrow\mathbb{R}$ be a continuous function. Then
\[
\frac{1}{n}\,\sup_{\theta\in\Theta}\Big|\sum_{i=1}^nh(X_i,\theta)\,e_i\Big|\ \asto\,0\,.
\]
\end{lemma}

\vspace*{1ex}\noindent
{\bf Proof.}\\
{\bf(a)} \ Denote \ $W_i\,:=\,|e_i|-{\rm E}\bigl(|e_i|\,\big|{\cal F}_{i-1}\bigr)$, $i\in\mathbb{N}$. It is easily seen that
the sequence of partial sums \ $\sum_{i=1}^n W_i$, $n\in\mathbb{N}$, is a martingale w.r.t. ${\cal F}_n$, $n\in\mathbb{N}$.
Since 
\[
{\rm E}\bigl(W_i^2\big|{\cal F}_{i-1}\bigr)\,=\,{\rm E}\bigl(e_i^2\big|{\cal F}_{i-1}\bigr)
\,-\,\bigl[{\rm E}\bigl(|e_i|\,\big|{\cal F}_{i-1}\bigr)\bigr]^2\,\le\,
{\rm E}\bigl(e_i^2\big|{\cal F}_{i-1}\bigr)\ \mbox{ a.s.}
\]
one has by (A2) \ $\sup_{i\in\mathbb{N}}{\rm E}\bigl(W_i^2\big|{\cal F}_{i-1}\bigr)\,<\infty\ \,\mbox{a.s.}$ and hence
\ $\sum_{i=1}^\infty i^{-2}{\rm E}\bigl(W_i^2\big|{\cal F}_{i-1}\bigr)\,<\infty\ \,\mbox{a.s.}$
By Theorem 2.18 of Hall and Heyde \cite{Hall-Heyde}, \ $\frac{1}{n}\sum_{i=1}^nW_i\,\asto 0$, i.e.,
\[
\frac{1}{n}\sum_{i=1}^n|e_i|\,-\,\frac{1}{n}\sum_{i=1}^n{\rm E}\bigl(|e_i|\,\big|{\cal F}_{i-1}\bigr)\ \asto\,0.
\]
By (A2) and Jensen's inequality \ $\sup_{i\in\mathbb{N}}{\rm E}\bigl(|e_i|\,\big|{\cal F}_{i-1}\bigr)\,<\infty \ \,\mbox{a.s.}$ \ 
from which one gets\\
$\limsup_{n\to\infty}\frac{1}{n}\sum_{i=1}^n|e_i|\,<\infty\ \,\mbox{a.s.}$\\[.5ex]
{\bf(b)} \ As it is easily seen, the sequence \ $\sum_{i=1}^nZ_ie_i$, $n\in\mathbb{N}$, is a martingale
 w.r.t. ${\cal F}_n$, $n\in\mathbb{N}$. By assumption there are two real random variables $U_1$ and $U_2$ such that 
$U_1\,=\,\sup_{i\in\mathbb{N}}{\rm E}\bigl(e_i^2\,\big|\,{\cal F}_{i-1}\bigr)$ a.s. and 
$U_2\,=\,\sup_{i\in\mathbb{N}}Z_i^2$ a.s. Hence 
\[
{\rm E}\bigl((Z_ie_i)^2\big|{\cal F}_{i-1}\bigr)\,\le\,U_2{\rm E}\bigl(e_i^2\big|{\cal F}_{i-1}\bigr)\,
\le\,U_2U_1 \ \mbox{ a.s. \ for all $i\in\mathbb{N}$.}
\]
So \ $\sum_{i=1}^\infty i^{-2}{\rm E}\bigl((Z_ie_i)^2\big|{\cal F}_{i-1}\bigr)\ <\infty\ \mbox{ a.s.}$
and the result follows from Theorem 2.18 of Hall and Heyde \cite{Hall-Heyde}.\\[.5ex]
{\bf(c)} \ Fix any $\alpha>0$.  By compactness of ${\cal X}\times\Theta$ and continuity of $h$ there exist a finite 
number $q\in\mathbb{N}$ and nonempty, pairwise disjoint, and measurable subsets $R_1,\ldots,R_q$ of ${\cal X}$ 
such that\  
$\bigcup_{j=1}^q R_j={\cal X}$ and  \ $|h(x,\theta)-h(z,\theta)|\le\alpha$ \ for all $x,z\in R_j$ and all $\theta\in\Theta$,
$1\le j\le q$. Choose any points $z_j^{(0)}\in R_j$, $1\le j\le q$, and denote $c_j(\theta)=h(z_j^{(0)},\theta)$, 
$1\le j\le q$, $\theta\in\Theta$. Then
\begin{equation}
c_j(\theta)-\alpha\,\le\,h(x,\theta)\,\le\,c_j(\theta)+\alpha\quad\forall\ x\in R_j,\ \forall\ \theta\in\Theta,\ 1\le j\le q.
\label{eq3-1}
\end{equation}
Introduce zero-one-valued random variables \ $Z_i^{(j)}:=\Ifkt\bigl(X_i^{-1}(R_j)\bigr)$, $i\in\mathbb{N}$, $1\le j\le q$,
i.e., $Z_i^{(j)}$ yields the value $1$ if the value of $X_i$ is in $R_j$, 
and otherwise $Z_i^{(j)}$ yields the value $0$ . 
Abbreviate $W_n(\theta):=\sum_{i=1}^nh(X_i,\theta)\,e_i$. Clearly, \ 
$W_n(\theta)=\sum_{j=1}^q\sum_{i=1}^nh(X_i,\theta)\,Z_i^{(j)}e_i$, and by (\ref{eq3-1}) for all $i$, $j$, and $\theta$,
\begin{equation}
\bigl(c_j(\theta)-{\rm sgn}(e_i)\,\alpha\bigr)\,Z_i^{(j)}e_i\,\le\,h(X_i,\theta)\,Z_i^{(j)}e_i\,\le\,
\bigl(c_j(\theta)+{\rm sgn}(e_i)\,\alpha\bigr)\,Z_i^{(j)}e_i \label{eq3-2}
\end{equation}
where ${\rm sgn}(t):=1$ if $t\ge0$ and ${\rm sgn}(t):=-1$ if $t<0$, for any real number $t$.  
Hence by summation in (\ref{eq3-2}) over $i$ and $j$,
\[
\sum_{j=1}^qc_j(\theta)\Bigl(\sum_{i=1}^nZ_i^{(j)}e_i\Bigr)\,-\alpha\sum_{i=1}^n|e_i|\,\le\,W_n(\theta)\,
\le\,\sum_{j=1}^qc_j(\theta)\Bigl(\sum_{i=1}^nZ_i^{(j)}e_i\Bigr)\,+\alpha\sum_{i=1}^n|e_i|.
\]
Denote \ $\overline{c}:=\sup_{(x,\theta)\in{\cal X}\times\Theta}|h(x,\theta)|$, which is finite and, clearly, 
$|c_j(\theta)|\le\overline{c}$ for all $j=1,\ldots,q$ and all $\theta\in\Theta$. Hence
\[
-\overline{c}\sum_{j=1}^q\Big|\sum_{i=1}^nZ_i^{(j)}e_i\Big|\,-\alpha\sum_{i=1}^n|e_i|\,\le\, W_n(\theta)
\,\le\,\overline{c}\sum_{j=1}^q\Big|\sum_{i=1}^nZ_i^{(j)}e_i\Big|\,+\alpha\sum_{i=1}^n|e_i|
\]
and thus
\[
\sup_{\theta\in\Theta}\big|W_n(\theta)\big|\,\le\,\overline{c}\sum_{j=1}^q\Big|\sum_{i=1}^nZ_i^{(j)}e_i\Big|\,+\alpha\sum_{i=1}^n|e_i|.
\]
Applying parts (a) and (b) of the lemma,
\[
\limsup_{n\to\infty}\Bigl(\frac{1}{n}\,\sup_{\theta\in\Theta}\big|W_n(\theta\big|\Bigr)\,\le\,\alpha\,U\quad\mbox{a.s.}
\]
where \ $U\,:=\,\limsup_{n\to\infty}\frac{1}{n}\sum_{i=1}^n|e_i|$,  which is almost surely finite. 
Since $\alpha>0$ was arbitrary the result follows.
\eop

\subsection{Adaptive GLM and  ML-estimators}
Now we specialize to an `adaptive generalized linear model' as follows. The parameter space
$\Theta$ is a compact subset of $\mathbb{R}^p$ provided with the usual Euclidean metric, the mean response function $\mu$ is of the form
\begin{equation}
\mu(x,\theta)\,=\,G\bigl(f^\trp(x)\,\theta\bigr), \ \ (x,\theta)\in{\cal X}\times\Theta, \label{eq3-3}
\end{equation}    
where $f\,:\,{\cal X}\longrightarrow\mathbb{R}^p$ is a given continuous function whose range $f({\cal X})$ spans $\mathbb{R}^p$
and $G\,:\,I\longrightarrow\mathbb{R}$ is a given continuously differentiable function on an open interval $I\subseteq\mathbb{R}$
with $\bigl\{f^\trp(x)\,\theta\,:\,(x,\theta)\in{\cal X}\times\Theta\bigr\}\subseteq I$ and whose derivative $G'$ is positive,
$G'(u)>0$ for all $u\in I$. 
The function $G$ is the inverse of the link function of the generalized linear model and $f^\trp(x)\,\theta$,
$(x,\theta)\in{\cal X}\times\Theta$, is the linear predictor.
Note that an interval may be unbounded from below or from above or both, where in the latter case the 
interval is the whole real line.  
Assumption (A2) is strengthened by an assumption (A2') below, stating 
that the conditional distribution of $Y_i$ given ${\cal F}_{i-1}$
belongs to a one-parameter exponential family of distributions $P_\tau$, $\tau\in J$, where $J\subseteq\mathbb{R}$ 
is an open interval. 
We employ the canonical (or `natural') parametrization of the one-parameter exponential family  
where  $\tau$ is its canonical parameter.
So $P_\tau$, $\tau\in J$, are probability distributions on the Borel sigma-field of the real line 
with densities w.r.t.~some Borel-measure $\nu$,
\begin{equation}
p_\tau(y)\,=\,K(y)\,\exp\bigl(\tau\,y\,-\,b(\tau)\bigr),\ \ y\in\mathbb{R},\ \ \tau\in J,\label{eq3-4}
\end{equation}
where $K$ is a nonnegative measurable function on $\mathbb{R}$ and $b$ is a real-valued function on $J$, which
is infinitely often differentiable, see e.g.~Fahrmeir and Kaufmann \cite{Fahrmeir-Kaufmann}, Section 2. 
In particular, the first and second derivatives of $b$ give the expectation and the variance of the distribution $P_\tau$,
resp., $b'(\tau)={\rm E}_{P_\tau}(Y)$ and $b''(\tau)={\rm Var}_{P_\tau}(Y)\,>0$ \ for all $\tau\in J$.
So the derivative $b'$ is a smooth and strictly increasing function and hence a bijection,  $b'\,:\,J\longrightarrow M$ where $M$ is the open 
interval of all expectations ${\rm E}_{P_\tau}(Y)$, $\tau\in J$. 
The inverse $(b')^{-1}$  assigns to each expectation $m\in M$ the 
parameter value $\tau=(b')^{-1}(m)\in J$ of the exponential family.    

Now, we assume the following (A2') which is stronger than assumption (A2) from Subsection 3.1. 
Recall that $\overline{\theta}$ denotes the fixed true parameter point.
\begin{itemize}
\item[{\bf(A2')}] The values of the inverse link function $G$ are contained in $M$, i.e.,
$G(I)\subseteq M$.\\ 
For each $i\in\mathbb{N}$ the conditional distribution of $Y_i$ given ${\cal F}_{i-1}$      
is equal to  $P_{\overline{\tau}_i}$ where\\ 
$\overline{\tau}_i=(b')^{-1}\bigl(G(f^\trp(X_i)\,\overline{\theta})\bigr)$. 
\end{itemize}
For the notion of a conditional distribution of a real-valued random variable given a sub-sigma-field 
we refer to \cite{Breiman}, p.~77, Definition 4.29. Note that 
$P_\tau$ has finite moments $m_k(\tau)={\rm E}_{P_\tau}\bigl(Y^k\bigr)$ of any order $k=1,2,\ldots$, and $m_k(\tau)$ 
is a continuous function of $\tau\in J$.  Assumption (A2') together with (A1) imply
the following. Firstly, 
${\rm E}\bigl(Y_i\big|\,{\cal F}_{i-1}\bigr)=m_1(\overline{\tau}_i)=  
G\bigl(f^\trp(X_i)\,\overline{\theta}\bigr)$. So, $e_i=Y_i -G\bigl(f^\trp(X_i)\,\overline{\theta}\bigr)$, $i\in\mathbb{N}$,
satisfy \ $Y_i=\mu(X_i,\overline{\theta})+e_i$ and ${\rm E}\bigl(e_i\big|\,{\cal F}_{i-1}\bigr)=0$, with $\mu$ from  
(\ref{eq3-3}). Secondly, ${\rm E}\bigl(e_i^2\big|\,{\cal F}_{i-1}\bigr)=m_2(\overline{\tau}_i)-\bigl(m_1(\overline{\tau}_i)\bigr)^2$
for all $i\in\mathbb{N}$, and since the values of all $\overline{\tau}_i$ are contained in some compact subinterval
of $J$ one has \ ${\rm E}\bigl(e_i^2\big|\,{\cal F}_{i-1}\bigr)\le C_2$ a.s. for all $i\in\mathbb{N}$ for some real constant $C_2>0$.
A similar conclusion holds for higher conditional moments of $e_i$, e.g. consider fourth moments:  
\[
{\rm E}\bigl(e_i^4\big|\,{\cal F}_{i-1}\bigr)=
\sum_{k=0}^4{4\choose k}\,m_k(\overline{\tau}_i)\,(-1)^{4-k}G^{4-k}\bigl(f^\trp(X_i)\,\overline{\theta}\bigr)\ \mbox{ a.s.},
\]
where $m_0(\overline{\tau}_i)=1$  and $G^0\bigl(f^\trp(X_i)\,\overline{\theta}\bigr)=1$.
The values of all the $X_i$, $i\in\mathbb{N}$, are in the  compact experimental region ${\cal X}$, and hence all the random variables 
$\overline{\tau}_i$, $i\in\mathbb{N}$, have their values in some compact subinterval of $J$. It follows that
\begin{equation}
{\rm E}\bigl(e_i^4\big|\,{\cal F}_{i-1}\bigr)\,\le C_4 \ \mbox{ a.s. } \ \mbox{for all $i\in\mathbb{N}$}\label{eq3-5}
\end{equation}
for some real constant $C_4>0$. 
To summarize: assumption (A2') together with (A1) imply (A2) and, moreover, (\ref{eq3-5}).
Obviously, this is due to the compactness of the experimental region ${\cal X}$
(and the continuity of $f$). Compactness of the parameter space $\Theta$, however, is not needed here since
(A2') as well as (A2) are local conditions at the true parameter point $\overline{\theta}$.

Fisher information matrices in a generalized linear model with univariate response whose observations follow
a one-parameter exponential family were derived in Atkinson and Woods \cite{Atkinson-Woods}, formula (13.3) on p.~473, and also for
the multivariate case in Fahrmeir and Kaufmann \cite{Fahrmeir-Kaufmann}, p.~347. Accordingly, 
we employ the following assumption (A3') on the family of functions 
$f_\theta$, $\theta\in\Theta$, defining the information matrices of designs via (\ref{eq1-0}).
\begin{itemize}
\item[{\bf(A3')}]
$f_\theta(x)\,=\,\ \varphi\bigl(f^\trp(x)\,\theta\bigr)\,f(x)$ \ for all $x\in{\cal X}$, $\theta\in\Theta$,\\
 where \ $\varphi(u)\,=\,G'(u)\Big/\sqrt{b''\Bigl((b')^{-1}\bigl(G(u)\bigr)\Bigr)}$, \ $u\in I$,
\end{itemize}
and where $f$ is a given continuous $\mathbb{R}^p$-valued function on ${\cal X}$ whose
range $f({\cal X})$ spans $\mathbb{R}^p$. 
In particular, by (A3') the family $f_\theta$, $\theta\in\Theta$, satisfies condition (GLM) from Section 2.
 
In what follows we focus on the asymptotics of adaptive maximum likelihood (ML) estimators. 
Note, however, that the adaptive estimators $\widehat{\theta}_n$, $n\ge n_{\rm\scriptsize st}$, employed by the algorithm 
may or may not be given by the adaptive ML-estimators $\widehat{\theta}_n^{(\rm\scriptsize ML)}$, $n\ge n_{\rm\scriptsize st}$.
The algorithm  may employ  any reasonable adaptive estimators $\widehat{\theta}_n$, $n\ge n_{\rm\scriptsize st}$,
e.g., the adaptive maximum quasi-likelihood estimators studied by Chen, Hu, and Ying \cite{Chen-Hu-Ying} in the case that the
function $G$ is defined on the whole real line, $I=\mathbb{R}$.
See also our remark below following Corollary \ref{cor3-1}. 
The main topics studied are strong consistency of the adaptive ML-estimators,
i.e., almost-sure convergence to the true parameter point $\overline{\theta}$, and asymptotic normality.
Strong consistency of the estimators $\widehat{\theta}_n$, $n\ge n_{\rm\scriptsize st}$, employed by the algorithm 
implies almost-sure asymptotic local D-optimality at $\overline{\theta}$ of the design sequence $\xi_n$
generated by the algorithm, which is an immediate consequence from Theorem \ref{theo2-2}. Note that the corollary
does not need any of the assumptions (A1), (A2), (A2'), or (A3').   
    
\begin{corollary}\quad\label{cor3-1}
If $\widehat{\theta}_n\asto\overline{\theta}$ then for any sequence $\widehat{\theta}_n'$ of estimators such that 
$\widehat{\theta}_n'\asto\overline{\theta}$ one has \ $M(\xi_n,\widehat{\theta}_n')\asto M(\xi^*_{\overline{\theta}},\overline{\theta})$,
where  $\xi^*_{\overline{\theta}}$ is a locally D-optimal design at $\overline{\theta}$.
\end{corollary}  

\noindent{\bf Remark.} \ 
Under assumptions (A1), (A2'), and (A3'), in the case $I=\mathbb{R}$ the adaptive maximum quasi-likelihood estimators   
studied by Chen, Hu, and Ying \cite{Chen-Hu-Ying}  turn out to be strongly consistent. 
In fact, by  Corollary \ref{cor2-1}, part (ii), and by (\ref{eq3-5}) one easily verifies the assumptions of Theorem 2 in \cite{Chen-Hu-Ying}
for the adaptive design sequence generated by the algorithm, irrespective of the employed sequence of adaptive estimators 
$\widehat{\theta}_n$, $n\ge n_{\rm\scriptsize st}$ in the algorithm. Our next result establishes strong consistency of
the adaptive maximum likelihood estimators, again irrespective of the employed sequence of estimators 
$\widehat{\theta}_n$, $n\ge n_{\rm\scriptsize st}$ in the algorithm.
\eop 

Assuming (A1) and (A2'), an adaptive ML-estimator 
$\widehat{\theta}_n^{(\rm\scriptsize ML)}=\widehat{\theta}_n^{(\rm\scriptsize ML)}(X_1,Y_1,\ldots,X_n,Y_n)$, 
for $n\ge n_{\rm\scriptsize st}$, is a maximizer of the log-likelihood
\begin{eqnarray} 
&&L_n(\theta)\,=\,\sum_{i=1}^n\Bigl(\log\bigl(K(Y_i)\bigr) + \tau_i(\theta)\,Y_i - b\bigl(\tau_i(\theta)\bigr)\Bigr), \label{eq3-8}\\
&&\mbox{ where }\ \tau_i(\theta)\,=\,(b')^{-1}\Bigl(G\bigl(f^\trp(X_i)\,\theta\bigr)\Bigr), \ \ 1\le i\le n.\label{eq3-8a}
\end{eqnarray} 
Note that with probability equal to one, $K(Y_i)>0$ for all $i\in\mathbb{N}$. Thus positivity of $K(Y_i)$, $i\in\mathbb{N}$, is assumed 
for the log-likelihood (\ref{eq3-8}). Note also that for the canonical link one gets $b'=G$ and hence
$\tau_i(\theta)=f^\trp(X_i)\,\theta$. 
The following result gives the strong consistency of the adaptive ML-estimators. 
  
\begin{theorem}\quad\label{theo3-1}
Under assumptions (A1), (A2'), and (A3'), one has \ $\widehat{\theta}_n^{(\rm\scriptsize ML)}\asto \overline{\theta}$. 
\end{theorem}

\noindent
{\bf Proof.} \    
For all $\theta\in\Theta$ one gets from (\ref{eq3-8}) and (\ref{eq3-8a}), observing that
$\tau_i(\overline{\theta})=\overline{\tau}_i$ and $Y_i=G\bigl(f^\trp(X_i)\,\overline{\theta}\bigr)+e_i$,
\begin{eqnarray}
&&L_n(\overline{\theta})-L_n(\theta)\,=\,\nonumber\\
&&\sum_{i=1}^n\bigl(\overline{\tau}_i-\tau_i(\theta)\bigr)\,G\bigl(f^\trp(X_i)\,\overline{\theta}\bigr)
+ \sum_{i=1}^n\bigl(\overline{\tau}_i-\tau_i(\theta)\bigr)e_i 
+ \sum_{i=1}^n\Bigl(b\bigl(\tau_i(\theta)\bigr)- b\bigl(\overline{\tau}_i\bigr)\Bigr).\label{eq3-9}
\end{eqnarray}
For each $i\in\mathbb{N}$, by second order Taylor expansion of $b(\tau)$ at $\overline{\tau}_i$,
\[
b\bigl(\tau_i(\theta)\bigr)\,=\,b\bigl(\overline{\tau}_i\bigr) + 
b'\bigl(\overline{\tau}_i\bigr)\,\bigl(\tau_i(\theta)-\overline{\tau}_i\bigr)
+{\textstyle\frac{1}{2}}\,b''\bigl(\widetilde{\tau}_i(\theta)\bigr)\,\bigl(\tau_i(\theta)-\overline{\tau}_i\bigr)^2
\]
with some $\widetilde{\tau}_i(\theta)$ from the interval whose end points are given by $\overline{\tau}_i$ and 
$\tau_i(\theta)$. Since 
\ $b'(\overline{\tau}_i)= G(f^\trp(X_i)\,\overline{\theta})$,  (\ref{eq3-9})
rewrites as
\begin{equation}
L_n(\overline{\theta})-L_n(\theta)\,=\,
\sum_{i=1}^n\bigl(\overline{\tau}_i-\tau_i(\theta)\bigr)e_i
+ \frac{1}{2}\sum_{i=1}^n b''(\widetilde{\tau}_i(\theta))\,\bigl(\tau_i(\theta)-\overline{\tau}_i\bigr)^2.
\label{eq3-10}
\end{equation}
By compactness of ${\cal X}$ and $\Theta$ and continuity of $f$ there is a compact subinterval $[c_1,c_2]\subseteq I$
such that $c_1\le f^\trp(x)\,\theta \le c_2$ for all $x\in{\cal X}$, $\theta\in\Theta$. Since $G$ and $(b')^{-1}$
are increasing functions,  \ $d_1\le (b')^{-1}\bigl(G(f^\trp(x)\,\theta)\bigr)\le d_2$ \ for all $x\in{\cal X}$, $\theta\in\Theta$, 
where $d_j=(b')^{-1}\bigl(G(c_j)\bigr)$, $j=1,2$, and $[d_1\,,\,d_2]\subseteq J$. In particular, 
 $d_1\le \tau_i(\theta)\le d_2$ for all $i$ and $\theta$. 
By continuity and positivity of $b''$
the minimum $\beta_0=\min_{d_1\le\tau\le d_2}b''(\tau)$ exists and $\beta_0>0$. Hence, in (\ref{eq3-10}), 
$b''\bigl(\widetilde{\tau}_i(\theta)\bigr)\ge\beta_0$ for all $i$ and $\theta$.
Since the composition $(b')^{-1}\bigl(G(u)\bigr)$, $u\in I$, is a continuously differentiable function
with positive derivative $H(u)\,=\,\frac{\rm d}{{\rm d}u}\bigl[ (b')^{-1}\bigl(G(u)\bigr)\bigr]$, $u\in I$, it follows that
\ $\beta_1\,=\,\min_{c_1\le u\le c_2} H(u)$ \
exists and $\beta_1>0$. By the mean value theorem  
\[
\big|(b')^{-1}\bigl(G(u_1)\bigr) - (b')^{-1}\bigl(G(u_2)\bigr)\big|\,\ge\, \beta_1\big|u_1-u_2\big|\quad\mbox{for all $u_1,u_2\in[c_1,c_2]$.}  
\] 
From (\ref{eq3-8a}) and  (\ref{eq3-10}) it follows that
\begin{equation}
L_n(\overline{\theta})-L_n(\theta)\,\ge\,\sum_{i=1}^n\bigl(\overline{\tau}_i-\tau_i(\theta)\bigr)e_i + 
{\textstyle\frac{1}{2}}\beta_0\beta_1^2\sum_{i=1}^n\bigl[f^\trp(X_i)\,\theta - f^\trp(X_i)\,\overline{\theta}\,\bigr]^2\ 
\mbox{ for all $\theta\in\Theta$.}\label{eq3-11}
\end{equation}
As in Wu \cite{Wu}, Lemma 1,
strong consistency of $\widehat{\theta}_n^{\rm\scriptsize (ML)}$ will follow if we prove that
for every $\delta>0$ such that the parameter subset 
$C(\overline{\theta},\delta)=\bigl\{\theta\in\Theta\,:\,\Vert\theta-\overline{\theta}\Vert\ge\delta\bigr\}$
is nonempty, one has
\[
\liminf_{n\to\infty}\Bigl(L_n(\overline{\theta})-\sup_{\theta\in C(\overline{\theta},\delta)} L_n(\theta)\Bigr)\,>0 \ \mbox{ a.s.}
\]
In fact, the $\liminf$ turns out to be equal to infinity almost surely, since we show that    
\begin{equation}
\liminf_{n\to\infty}\frac{1}{n}\Bigl(L_n(\overline{\theta})-\sup_{\theta\in C(\overline{\theta},\delta)} L_n(\theta)\Bigr)\ >0
\ \mbox{a.s.}\label{eq3-12}
\end{equation}
By (\ref{eq3-11}) and the trivial inequality $a\ge-|a|$ for all real $a$,
\begin{eqnarray}
&&\frac{1}{n}\Bigl(L_n(\overline{\theta})-\sup_{\theta\in C(\overline{\theta},\delta)} L_n(\theta)\Bigr)\,\ge\nonumber\\
&&-\frac{1}{n}\sup_{\theta\in\Theta}\Big|\sum_{i=1}^n\bigl(\overline{\tau}_i-\tau_i(\theta)\bigr)e_i\Big|
+ {\textstyle\frac{1}{2}}\beta_0\beta_1^2
\frac{1}{n}\inf_{\theta\in C(\overline{\theta},\delta)}\sum_{i=1}^n\bigl[f^\trp(X_i)\,\theta - f^\trp(X_i)\,\overline{\theta}\,\bigr]^2
\label{eq3-13}
\end{eqnarray}
Since \ $\overline{\tau}_i-\tau_i(\theta)=h(X_i,\theta)$ for all $i\in\mathbb{N}$ and $\theta\in\Theta$,
with $h(x,\theta)=(b')^{-1}\bigl(f^\trp(x)\,\overline{\theta}\bigr) - (b')^{-1}\bigl(f^\trp(x)\,\theta\bigr)$,
$(x,\theta)\in{\cal X}\times\Theta$, it follows by Lemma \ref{lem3-1}, part (c), that
\begin{equation}
\frac{1}{n}\sup_{\theta\in\Theta}\Big|\sum_{i=1}^n\bigl(\overline{\tau}_i-\tau_i(\theta)\bigr)e_i\Big|\ \asto 0.
\label{eq3-13a}
\end{equation}
It remains to show that the $\liminf$ of the second term on the r.h.s. of (\ref{eq3-13}) is positive almost surely.
Consider an arbitrary path of the adaptive process and, in particular, a path $x_i$, $i\in \mathbb{N}$, of
the sequence of random variables $X_i$, $i\in\mathbb{N}$. 
With the generated design sequence
$\xi_n$, $n\ge n_{\rm\scriptsize st}$, we can write, for all $n\ge n_{\rm\scriptsize st}$, 
\begin{equation}
\frac{1}{n}\inf_{\theta\in C(\overline{\theta},\delta)}\sum_{i=1}^n\bigl[f^\trp(x_i)\,\theta - f^\trp(x_i)\,\overline{\theta}\,\bigr]^2
\,=\,\inf_{\theta\in C(\overline{\theta},\delta)}\int_{\cal X} \bigl[f^\trp(x)\,(\theta-\overline{\theta})\bigr]^2\,{\rm d}\xi_n(x).
\label{eq3-13b}
\end{equation}
For any $\theta\in\Theta$, $\theta\not=\overline{\theta}$, denote $c_\theta=(\theta-\overline{\theta})/\Vert\theta-\overline{\theta}\Vert$
and  $V_{p-1,\theta}=\bigl\{a\in\mathbb{R}^p\,:\,c^\trp_\theta a=0\bigr\}$. By Theorem \ref{theo2-1} there exist $n_0\ge n_{\rm\scriptsize st}$,
$\varepsilon>0$, and $\alpha\in(\,0\,,\,1\,)$ such that 
\[
\xi_n\Bigl(f_{\theta_n}^{-1}\bigl(\overline{V}_{p-1,\theta}(\varepsilon)\bigr)\Bigr)\,\le\,\alpha\ \mbox{ for all $\theta\not=\overline{\theta}$ and
$n\ge n_0$.}
\]
Using (A3') and $f^\trp(x)\,\theta\in[\,c_1\,,\,c_2\,]\subseteq I$ for all $(x,\theta)\in{\cal X}\times\Theta$, let   
\[
\varphi_{\rm\scriptsize min}:=\inf_{c_1\le u\le c_2}\varphi(u)\quad\mbox{and}\quad \varphi_{\rm\scriptsize max}:=\sup_{c_1\le u\le c_2}\varphi(u),
\]
hence $0<\varphi_{\rm\scriptsize min}\le \varphi_{\rm\scriptsize max}<\infty$. 
Define \ $\varepsilon':=\varepsilon\,\varphi_{\rm\scriptsize min}/\varphi_{\rm\scriptsize max}$. As in the proof of Lemma \ref{lem2-5}, part (i),
one gets  
\[
f^{-1}\bigl(\overline{V}_{p-1,\theta}(\varepsilon')\bigr)\subseteq f_{\theta_n}^{-1}\bigl(\overline{V}_{p-1,\theta}(\varepsilon)\bigr)
\ \mbox{ for all $n\ge n_{\rm\scriptsize st}$ and all $\theta\not=\overline{\theta}$.}
\]
Note that \ $f^{-1}\bigl(\overline{V}_{p-1,\theta}(\varepsilon')\bigr)\,=\,\bigl\{x\in{\cal X}\,:\,|c^\trp_\theta f(x)|\le
\varepsilon'\bigr\}$.
\ Taking the complementary sets and observing that $|c^\trp_\theta f(x)|>\varepsilon'$ is equivalent to
$\big|f^\trp(x)\,(\theta-\overline{\theta})\big|> \varepsilon'\Vert \theta-\overline{\theta}\Vert$, which in the case 
$\theta\in C(\overline{\theta},\delta)$ implies $\big|f^\trp(x)\,(\theta-\overline{\theta})\big|> \varepsilon'\delta$,
we have
\[
\xi_n\Bigl(\bigl\{x\in{\cal X}\,:\,
\big|f^\trp(x)\,(\theta-\overline{\theta})\big|> \varepsilon'\delta\bigr\}\Bigr) \,\ge\,1-\alpha\ 
\mbox{ for all $\theta\in C(\overline{\theta},\delta)$ and $n\ge n_0$.}
\]
Hence
\[
\inf_{\theta\in C(\overline{\theta},\delta)}\int_{\cal X} \bigl[f^\trp(x)\,(\theta-\overline{\theta})\bigr]^2\,{\rm d}\xi_n(x)\,\ge\,
(\varepsilon'\delta)^2(1-\alpha)\,>0\ \mbox{ for all $n\ge n_0$.}
\]
Together with (\ref{eq3-13b}), (\ref{eq3-13a}), and (\ref{eq3-13}) the proof of (\ref{eq3-12}) is complete and
$\widehat{\theta}_n^{(\rm\scriptsize ML)}\asto \overline{\theta}$ follows.   
\eop  

The next result shows the asymptotic normality of adaptive ML-estimators if the true parameter point $\overline{\theta}$
is an interior point of the parameter space $\Theta$, i.e., there exists a $\rho>0$ such that
\begin{equation}
B(\overline{\theta},\rho):=\bigl\{a\in\mathbb{R}^p\,:\,\Vert a-\overline{\theta}\Vert <\rho\bigr\}\subseteq\Theta.
\label{eq3-14aa}
\end{equation}
The $p$-dimensional normal distribution with expectation $0$ and (positive definite) covariance matrix $C$ is denoted by 
${\rm N}(0,C)$. For $C=I_p$, the $p\times p$ identity matrix, ${\rm N}(0,I_p)$ is  the $p$-dimensional standard normal distribution.  
For a sequence $W_n$ of $\mathbb{R}^p$-valued random variables, 
convergence in distribution of $W_n$ (as $n\to \infty$) to a $p$-dimensional normal distribution ${\rm N}(0,C)$ 
is abbreviated by $W_n\dto{\rm N}(0,C)$.   
In the next theorem, the assumption of strong consistency of the adaptive estimators $\widehat{\theta}_n$ employed by the algorithm
is met, by Theorem \ref{theo3-1}, if $\widehat{\theta}_n=\widehat{\theta}_n^{\rm\scriptsize (ML)}$,    
$n\ge n_{\rm\scriptsize st}$. 
 
\begin{theorem}\quad \label{theo3-2}  
Assume (A1), (A2'), and (A3'). 
Assume further that  $G$ is twice continuously differentiable, $\overline{\theta}$ is an interior point of $\Theta$, and 
$\widehat{\theta}_n\asto\overline{\theta}$. Then:
\[
\sqrt{n}\,M^{1/2}\bigl(\xi_n,\widehat{\theta}_n^{\rm\scriptsize (ML)}\bigr)\,\bigl(\widehat{\theta}_n^{\rm\scriptsize (ML)}-\overline{\theta}\bigr)
\,\dto\,{\rm N}(0,I_p).
\]
Also, denoting by $M_*=M\bigl(\xi^*_{\overline{\theta}},\overline{\theta})$ the information matrix of a locally D-optimal  
design at $\overline{\theta}$, one has
\[
\sqrt{n}\,\bigl(\widehat{\theta}_n^{\rm\scriptsize (ML)}-\overline{\theta}\bigr)
\,\dto\,{\rm N}\bigl(0,M_*^{-1}\bigr).
\]
\end{theorem}

\noindent
{\bf Proof.} \ Choose a positive $\rho_1<\rho$ where  $\rho>0$ is according to (\ref{eq3-14aa}). Then the compact ball
$\overline{B}(\overline{\theta},\rho_1)=\bigl\{a\in\mathbb{R}^p\,:\,\Vert a-\overline{\theta}\Vert\le\rho_1\bigr\}$
is contained in the interior of $\Theta$. By Theorem \ref{theo3-1} 
$\widehat{\theta}_n^{\rm\scriptsize (ML)}\in\overline{B}(\overline{\theta},\rho_1)\ \mbox{ a.s.}$ \ 
if $n$ is large enough, i.e.,
if $n$ is greater than or equal to the value of some random variable $N$ whose values are in $\mathbb{N}$. 
Let $n\ge n_{\rm\scriptsize st}$ be given and
denote by  $\{N\le n\}$ the event (subset of $\Omega$) that the random variable $N$
yields a value less than or equal to $n$.  
Consider the log-likelihood function $L_n(\theta)$ from (\ref{eq3-8}), (\ref{eq3-8a}) and its gradient 
$S_n(\theta)=\nabla L_n(\theta)$ w.r.t. $\theta$,
which is often called the score function, for $\theta\in\overline{B}(\overline{\theta},\rho_1)$. One obtains  
\begin{eqnarray}
&&S_n(\theta)\,=\,\sum_{i=1}^n\Bigl(Y_i-G\bigl(f^\trp(X_i)\,\theta\bigr)\Bigr)\,H\bigl(f^\trp(X_i)\,\theta\bigr)\,f(X_i),
\label{eq3-14}\\
&&\mbox{ where }\ H(u)\,=\,\frac{G'(u)}{b''\Bigl((b')^{-1}\bigl(G(u)\bigr)\Bigr)}\ \mbox{ for all $u\in I$.}
\label{eq3-14a}
\end{eqnarray}
Abbreviate
\[
R_i(\theta)\,=\,\Bigl(Y_i-G\bigl(f^\trp(X_i)\theta\bigr)\Bigr)\,H\bigl(f^\trp(X_i)\theta\bigr), \ (1\le i\le n).
\]
Since $S_n\bigl(\widehat{\theta}_n^{\rm\scriptsize (ML)}\bigr)=0$ \ a.s. on $\{N\le n\}$ \ we get from (\ref{eq3-14}),
\begin{equation}
-S_n(\overline{\theta})=S_n\bigl(\widehat{\theta}_n^{\rm\scriptsize (ML)}\bigr)-S_n(\overline{\theta})=
\sum_{i=1}^n \Bigl(R_i\bigl(\widehat{\theta}_n^{\rm\scriptsize (ML)}\bigr)-R_i(\overline{\theta})\Bigr)\,f(X_i)
\ \mbox{ a.s. on $\{N\le n\}$.}
\label{eq3-15}
\end{equation}
The function $H$ from (\ref{eq3-14a}) is continuously differentiable. Denote its derivative by $H'$.
The gradient (w.r.t. $\theta$) of $R_i(\theta)$ is given by
\begin{equation}
\nabla R_i(\theta)\,=\,\Bigl(-G'\bigl(f^\trp(X_i)\,\theta\bigr)\,H\bigl(f^\trp(X_i)\,\theta\bigr) +
\bigl[Y_i-G\bigl(f^\trp(X_i)\,\theta\bigr]\,H'\bigl(f^\trp(X_i)\theta\bigr)\Bigr)\,f(X_i).
\label{eq3-15a}
\end{equation}
By the mean value theorem, for each $i=1,\ldots,n$  there is some $\widetilde{\theta}_{i,n}$ on the line segment joining 
$\overline{\theta}$ and $\widehat{\theta}_n^{\rm\scriptsize (ML)}$ such that
\[
R_i\bigl(\widehat{\theta}_n^{\rm\scriptsize (ML)}\bigr)-R_i(\overline{\theta})\,=\,
\bigl(\nabla R_i\bigr)^\trp(\widetilde{\theta}_{i,n})\,\bigl(\widehat{\theta}_n^{\rm\scriptsize (ML)}-\overline{\theta}\bigr).
\]
Together with (\ref{eq3-15}) and (\ref{eq3-15a}), and observing that \ $G'\,H\,=\,\varphi^2$, we get
\begin{eqnarray}
 -S_n(\overline{\theta}) &=&\bigl(-A_n+B_n\bigr)\,
\bigl(\widehat{\theta}_n^{\rm\scriptsize (ML)}-\overline{\theta}\bigr)
\ \mbox{ a.s. on $\{N\le n\}$,}
\ \mbox{ where}\label{eq3-16}\\
A_n &=& \sum_{i=1}^n\varphi^2\bigl(f^\trp(X_i)\,\widetilde{\theta}_{i,n}\bigr)\,
f(X_i)\,f^\trp(X_i)\label{eq3-16a},\\
B_n &=& \sum_{i=1}^n\Bigl(Y_i-G\bigl(f^\trp(X_i)\,\widetilde{\theta}_{i,n}\bigr)\Bigr)\,
H'\bigl(f^\trp(X_i)\,\widetilde{\theta}_{i,n}\bigr)\,f(X_i)\,f^\trp(X_i). \label{eq3-16b}
\end{eqnarray}
Since $M(\xi_n,\theta)=\frac{1}{n}\sum_{i=1}^n\varphi^2\bigl(f^\trp(X_i)\,\theta\bigr)\,f(X_i)\,f^\trp(X_i)$ 
we can write
\begin{eqnarray}
A_n &=& n\,M(\xi_n,\overline{\theta})\,+\,D_n,\ \mbox{ where }\label{eq3-17}\\
D_n &=& \sum_{i=1}^n\bigl[\varphi^2\bigl(f^\trp(X_i)\,\widetilde{\theta}_{i,n}\bigr)-
\varphi^2\bigl(f^\trp(X_i)\,\overline{\theta}\bigr)\bigr]\,f(X_i)\,f^\trp(X_i).\label{eq3-17a}
\end{eqnarray}
So, by (\ref{eq3-16}) after some slight manipulations, 
\begin{equation}
\frac{1}{\sqrt{n}}\,M_*^{-1/2}\,S_n(\overline{\theta})\,=\,
M_*^{-1/2}\,\Bigl[M(\xi_n,\overline{\theta})+
\frac{1}{n}D_n - \frac{1}{n}B_n\Bigr]\,\bigl[\sqrt{n}\bigl(\widehat{\theta}_n^{\rm\scriptsize (ML)}-\overline{\theta}\bigr)\bigr]
\mbox{ a.s. on $\{N\le n\}$.}
\label{eq3-18}
\end{equation}
Next we show that
\begin{equation}
\frac{1}{\sqrt{n}}\,M_*^{-1/2}\,S_n(\overline{\theta})\,\dto\,{\rm N}(0,I_p).
\label{eq3-19}
\end{equation}
Regarding the Cram\'{e}r-Wold device let any $v\in\mathbb{R}^p$ with $\Vert v\Vert=1$ be given. 
Using (\ref{eq3-14}) for $\theta=\overline{\theta}$ and inserting $Y_i=G\bigl(f^\trp(X_i)\,\overline{\theta}\bigr)\,+e_i$
we can write
\begin{eqnarray*}
&&\frac{1}{\sqrt{n}}\,v^\trp M_*^{-1/2}\,S_n(\overline{\theta})\,=\,
\frac{1}{\sqrt{n}}\,\sum_{i=1}^n\widetilde{e}_i,\ \mbox{ where }\\
&&\widetilde{e}_i\,=\,e_iH\bigl(f^\trp(X_i)\,\overline{\theta}\bigr)\,v^\trp M_*^{-1/2}\,f(X_i).
\end{eqnarray*}
Clearly, for each $i\in\mathbb{N}$ the random variable 
$Z_i:=H\bigl(f^\trp(X_i)\,\overline{\theta}\bigr)\,v^\trp M_*^{-1/2}f(X_i)$ is 
${\cal F}_{i-1}$-measurable and $|Z_i|\le c$ for all $i\in\mathbb{N}$ for some finite constant $\widetilde{c}$. 
Since  $\widetilde{e}_i=e_iZ_i$, $i\in\mathbb{N}$, one easily verifies that    
the sequence of partial sums $\sum_{i=1}\widetilde{e}_i$, $n\in\mathbb{N}$, is a martingale w.r.t. ${\cal F}_n$,
$n\in\mathbb{N}$. By Corollary 3.1 (p.~58) in \cite{Hall-Heyde} the convergence \  
$\frac{1}{\sqrt{n}}\,\sum_{i=1}^n\widetilde{e}_i\dto {\rm N}(0,1)$ \ holds if the following two conditions (A) and (B) are satisfied.
\begin{itemize}
\item[(A)] $\displaystyle\frac{1}{n}\sum_{i=1}^n{\rm E}\bigl(\widetilde{e}_i^2\big|\,{\cal F}_{i-1}\bigr)\,\asto 1$.
\item[(B)] 
$\displaystyle\frac{1}{n}\sum_{i=1}^n{\rm E}\Bigl(\widetilde{e}_i^2\Ifkt\bigl(|\widetilde{e}_i|>\sqrt{n}\,\varepsilon\bigr)\big|\,{\cal F}_{i-1}\Bigr)
\,\asto 0$ \ for all $\varepsilon>0$. 
\end{itemize}
\underline{Ad (A)}. \ ${\rm E}\bigl(\widetilde{e}_i^2\big|\,{\cal F}_{i-1}\bigr)=
{\rm E}\bigl(e_i^2\big|\,{\cal F}_{i-1}\bigr)\,Z_i^2$ and
$Z_i^2\,=\,H^2\bigl(f^\trp(X_i)\,\overline{\theta}\bigr)\,
v^\trp M_*^{-1/2}\,f(X_i)\,f^\trp(X_i)\,M_*^{-1/2}\,v$.
By assumption (A2'), 
${\rm E}\bigl(e_i^2\big|\,{\cal F}_{i-1}\bigr)=b''\Bigl((b')^{-1}\bigl(G(f^\trp(X_i)\,\overline{\theta})\bigr)\Bigr)$.
Together with (\ref{eq3-14a}) this yields \ 
${\rm E}\bigl(\widetilde{e}_i^2\big|\,{\cal F}_{i-1}\bigr)=\varphi^2\bigl(f^\trp(X_i)\,\overline{\theta}\bigr)\,  
v^\trp M_*^{-1/2}\,f(X_i)\,f^\trp(X_i)\,M_*^{-1/2}v$,
and hence
\[
\frac{1}{n}\sum_{i=1}^n{\rm E}\bigl(\widetilde{e}_i^2\big|\,{\cal F}_{i-1}\bigr)\,=\,
v^\trp M_*^{-1/2}\,
M(\xi_n,\overline{\theta})\,M_*^{-1/2}v\,\asto 1,
\]
where the convergence follows from Corollary \ref{cor3-1}.\\[.5ex]
\underline{Ad (B)}.\ Using the trivial inequality \
$\widetilde{e}_i^2\,\Ifkt\bigl(|\widetilde{e}_i|>\sqrt{n}\,\varepsilon\bigr)\,\le\,\frac{1}{\varepsilon^2 n}\,\widetilde{e}_i^4$ \
we obtain
\begin{eqnarray*}
&&\frac{1}{n}\sum_{i=1}^n{\rm E}\Bigl(\widetilde{e}_i^2\Ifkt\bigl(|\widetilde{e}_i|>\sqrt{n}\,\varepsilon\bigr)\big|\,{\cal F}_{i-1}\Bigr)\,
\le\,\frac{1}{\varepsilon^2n^2}\,\sum_{i=1}^n{\rm E}\bigl(\widetilde{e}_i^4\big|\,{\cal F}_{i-1}\bigr)\\
&&=\, \frac{1}{\varepsilon^2n^2}\,\sum_{i=1}^nZ_i^4{\rm E}\bigl(e_i^4\big|\,{\cal F}_{i-1}\bigr)\,\le\,
\frac{c^4\,C_4}{\varepsilon^2n}\,\asto 0,
\end{eqnarray*}
where we have used (\ref{eq3-5}). \\[.5ex]
So by the Cram\'{e}r-Wold device (\ref{eq3-19}) follows and hence, by (\ref{eq3-18}),
\begin{equation}
M_*^{-1/2}\,\Bigl[M(\xi_n,\overline{\theta})+
\frac{1}{n}D_n - \frac{1}{n}B_n\Bigr]\,\bigl[\sqrt{n}\bigl(\widehat{\theta}_n^{\rm\scriptsize (ML)}-\overline{\theta}\bigr)\bigr]
\,\dto{\rm N}(0,I_p). 
\label{eq3-20}
\end{equation}
Next we show that 
\begin{equation}
\frac{1}{n}D_n\asto 0\quad\mbox{and}\quad \frac{1}{n}B_n\asto 0. \label{eq3-21}
\end{equation}
By (\ref{eq3-17a}), 
\begin{eqnarray*}
&&\frac{1}{n}D_n\,=\,\frac{1}{n}\sum_{i=1}^n\bigl[\varphi^2\bigl(f^\trp(X_i)\,\widetilde{\theta}_{i,n})\bigr) -
\varphi^2\bigl(f^\trp(X_i)\,\overline{\theta})\bigr)\bigr]\,f(X_i)\,f^\trp(X_i),\\
&&\mbox{hence }\ 
\Big\Vert\frac{1}{n}D_n\Big\Vert\,\le\,\max_{1\le i\le n}\Big|  
\varphi^2\bigl(f^\trp(X_i)\,\widetilde{\theta}_{i,n})\bigr) -
\varphi^2\bigl(f^\trp(X_i)\,\overline{\theta}\bigr)\big|\ \frac{1}{n}\sum_{i=1}^n\big\Vert f(X_i)\,f^\trp(X_i)\big\Vert,
\end{eqnarray*}
Note that \ $\big\Vert f(X_i)\,f^\trp(X_i)\big\Vert=\big\Vert f(X_i)\big\Vert^2\le\gamma_0^2$ for all $i\in\mathbb{N}$, 
where $\gamma_0=\sup_{x\in{\cal X}}\Vert f(x)\Vert<\infty$. 
By  $f^\trp(X_i)\,\theta\in[\,c_1\,,\,c_2\,]\subseteq I$ for all $i\in\mathbb{N}$ and
$\theta\in\overline{B}(\overline{\theta},\rho_1)$ for some compact interval $[\,c_1\,,\,c_2\,]$, 
by the uniform continuity of the function $\varphi^2$ on the compact interval, and by
\[
\max_{1\le i\le n}\big|f^\trp(X_i)\,\widetilde{\theta}_{i,n}-f^\trp(X_i)\,\overline{\theta}\big|\,
\le \gamma_0\Vert\widetilde{\theta}_{i,n}-\overline{\theta}\Vert\le 
\gamma_0\Vert\widehat{\theta}_n^{\rm\scriptsize(ML)}-\overline{\theta}\Vert\asto 0
\]
it follows that
\[
\max_{1\le i\le n}\big|\varphi^2\bigl(f^\trp(X_i)\,\widetilde{\theta}_{i,n}\bigr) -
\varphi^2\bigl(f^\trp(X_i)\,\overline{\theta}\bigr)\big|\,\asto 0.
\]
From this the first convergence statement in (\ref{eq3-21}) follows. To prove the second convergence statement 
in (\ref{eq3-21}) we write $Y_i-G\bigl(f^\trp(X_i)\,\widetilde{\theta}_{i,n}\bigr)=
G\bigl(f^\trp(X_i)\,\overline{\theta}\bigr) -G\bigl(f^\trp(X_i)\,\widetilde{\theta}_{i,n}\bigr) + e_i$, and  
together with the definition  of $B_n$ in (\ref{eq3-16b}), 
\begin{eqnarray*}
B_n &=&\,B_n^{(1)} + B_n^{(2)} + B_n^{(3)},\ \mbox{ where}\\
B_n^{(1)} &=& \sum_{i=1}^n\bigl[G\bigl(f^\trp(X_i)\,\overline{\theta}\bigr) -G\bigl(f^\trp(X_i)\,\widetilde{\theta}_{i,n}\bigr)\bigr]\,
H'\bigl(f^\trp(X_i)\,\widetilde{\theta}_{i,n}\bigr)\,f(X_i)\,f^\trp(X_i),\\
B_n^{(2)} &=& \sum_{i=1}^n e_i\,H'\bigl(f^\trp(X_i)\,\overline{\theta}\bigr)\,f(X_i)\,f^\trp(X_i),\\
B_n^{(3)} &=& \sum_{i=1}^n e_i\,
\bigl[H'\bigl(f^\trp(X_i)\,\widetilde{\theta}_{i,n}\bigr)-H'\bigl(f^\trp(X_i)\,\overline{\theta}\bigr)\bigr]\,
f(X_i)\,f^\trp(X_i).
\end{eqnarray*}
One concludes \ $\frac{1}{n}B_n^{(1)}\asto 0$ \ by similar arguments as used above when showing $\frac{1}{n}D_n\asto 0$ where, in particular, 
uniform continuity of $G$ on the compact interval $[\,c_1\,,\,c_2\,]\subseteq I$ and boundedness of $H'$ on that interval  
are utilized.  
Consider the sequence of matrices $\frac{1}{n}B_n^{(2)}$ entrywise. The $(k,\ell)$-th entry (where $1\le k,\ell\le p$)
of $\frac{1}{n}B_n^{(2)}$ has the form \ $\frac{1}{n}\sum_{i=1}^n Z_i\,e_i$ where \ 
$Z_i=H'\bigl(f^\trp(X_i)\,\overline{\theta}\bigr)\,f_k(X_i)\,f_\ell(X_i)$. Note that $|Z_i|\le c$ 
for all $i\in\mathbb{N}$ for some real constant $c>0$. 
From Lemma \ref{lem3-1}, part (b), it follows that \ $\frac{1}{n}\sum_{i=1}^n Z_i\,e_i\,\asto 0$. Hence \ 
$\frac{1}{n}B_n^{(2)}\asto 0$. For $\frac{1}{n}B_n^{(3)}$ we consider again the $(k,\ell)$-th entry for any given $1\le k,\ell\le p$,
and we obtain
\begin{eqnarray*}
&&\Big|\frac{1}{n}\sum_{i=1}^n e_i\,\bigl[H'\bigl(f^\trp(X_i)\,\widetilde{\theta}_{i,n}\bigr)-H'\bigl(f^\trp(X_i)\,\overline{\theta}\bigr)\bigr]
\,f_k(X_i)f_\ell(X_i)\Big|\\
&&\le\,\gamma_0^2\,\max_{1\le i\le n}\big|H'\bigl(f^\trp(X_i)\,\widetilde{\theta}_{i,n}\bigr)-H'\bigl(f^\trp(X_i)\,\overline{\theta}\bigr)\big|\,
\frac{1}{n}\sum_{i=1}^n|e_i|.
\end{eqnarray*}
The uniform continuity of $H'$ on $[\,c_1\,,\,c_2\,]$ yields \
$\max_{1\le i\le n}\big|H'\bigl(f^\trp(X_i)\,\widetilde{\theta}_{i,n}\bigr)-H'\bigl(f^\trp(X_i)\,\overline{\theta}\bigr)\big|\asto 0$.
Together with Lemma \ref{lem3-1}, part (a), it follows that the $(k,\ell)$-th entry of $\frac{1}{n}B_n^{(3)}$ converges to zero almost
surely and hence $\frac{1}{n}B_n^{(3)}\asto 0$. So we have proved (\ref{eq3-21}).  
By Corollary \ref{cor3-1}, $M(\xi_n,\overline{\theta})\asto M_*$ and hence   
\[
M_*^{-1/2}\,\Bigl[M(\xi_n,\overline{\theta})+
\frac{1}{n}D_n - \frac{1}{n}B_n\Bigr]\,\asto M_*^{1/2}.
\]
Together with (\ref{eq3-20}), observing  $\lim_{n\to\infty}\mathbb{P}_{\overline{\theta}}\bigl(\{N\le n\}\bigr)=1$
and using standard properties of convergence in distribution, it follows that 
for any sequence of $p\times p$ random matrices $Q_n$ such that $Q_n\asto M_*^{1/2}$
one has \ $Q_n\bigl[\sqrt{n}\bigl(\widehat{\theta}_n^{\rm\scriptsize (ML)}-\overline{\theta}\bigr)\bigr]\dto {\rm N}(0,I_p)$.
In particular, the convergence holds for the sequence $Q_n=M^{1/2}(\xi_n,\widehat{\theta}_n^{\rm\scriptsize (ML)})$
and the constant sequence $Q_n=M_*^{1/2}$. Hence the result follows. 
\eop

\begin{appendix}
\setcounter{equation}{0}
\section{Appendix: known auxiliary results}
Four well-known results on nonnegative definite matrices are stated below, 
which have been used throughout the proofs.
If $A$ is a (real) $p\times q$ matrix then the range of $A$ is given by ${\rm range}(A)=\{Az\,:\,z\in\mathbb{R}^q\}$.  
A generalized inverse of a $p\times q$ matrix $A$ is denoted by
$A^-$ which is by definition a $q\times p$ matrix satisfying $AA^-A=A$. As it is easily seen, 
if $A$ is symmetric $p\times p$ and $b\in{\rm range}(A)$ 
then the value $b^\trp A^-b$ is the same for all generalized inverses of $A$. 

\begin{itemize}
\item[(M1)]
If $A$ and $B$ are nonnegative definite $p\times p$ matrices and $A\le B$, then  \ ${\rm range}(A)\subseteq{\rm range}(B)$ \ 
and \ $z^\trp A^- z\,\ge z^\trp B^- z$ \ for all $z\in{\rm range}(A)$. 
See Stepniak, Wang and Wu \cite{Stepniak-Wang-Wu}, Lemma 2. 
\item[(M2)]
If $a_1,\ldots,a_r\in\mathbb{R}^p$ and $\lambda_1,\ldots,\lambda_r\in(\,0\,,\,\infty)$ such that $\sum_{j=1}^r\lambda_j=1$, then 
\[
\Bigl(\sum_{j=1}^r\lambda_ja_j\Bigr)^\trp \Bigl(\sum_{j=1}^r\lambda_j\,a_ja_j^\trp\Bigr)^- 
\Bigl(\sum_{j=1}^r\lambda_ja_j\Bigr)
\,\le\,\sum_{j=1}^r \lambda_j\,a_j^\trp\bigl(a_ja_j^\trp\bigr)^-a_j\,\le\,1.
\]
The first inequality is a special case of Theorem 4.2 in Gaffke and Krafft \cite{Gaffke-Krafft}; the second inequality
follows from $a^\trp(aa^\trp)^-a=1$ if $a\not=0$, and  $a^\trp(aa^\trp)^-a=0$ if $a=0$, for any $a\in\mathbb{R}^p$.
\item[(M3)]
If $M$ is a positive definite $p\times p$ matrix and $v\in\mathbb{R}^p$, $v\not=0$, then
\[
v^\trp M^{-1}v\,=\,\max\Bigl\{\frac{1}{b^\trp Mb}\,:\,b\in\mathbb{R}^p,\ v^\trp b=1\Bigr\}\ ,
\]
and the maximum on the right hand side is attained for \ $b_0\,=\,M^{-1}v/(v^\trp M^{-1}v)$.
Actually, the inequality is a special case of a more general matrix inequality, 
see Pukelsheim \cite{Puk}, Section 1.21.  
\item[(M4)]
Let $A$ be  a $p\times q$ matrix with columns  $a_1,\ldots,a_q\in\mathbb{R}^p$. Then 
\[
\det(A^\trp A)\,=\,\prod_{j=1}^q {\rm dist}^2(a_j,V_{j-1}),\ \mbox{ where } V_0:=\{0\},\ 
V_k:={\rm span}\{a_1,\ldots,a_k\},\ 1\le k\le q-1,
\]
and where ${\rm dist}^2(a,V)=\inf_{v\in V}\Vert a-v\Vert^2$ denotes the squared Euclidean distance of a vector $a\in\mathbb{R}^p$ and
a linear subspace $V$ of $\mathbb{R}^p$. 
In fact, the formula trivially holds if the vectors $a_1,\ldots,a_q$ are linearly dependent, in
which case both sides of the formula are equal to zero. Also, the case $q=1$ is trivial. 
Let $q\ge2$ and let $a_1,\ldots,a_q$ be linearly independent. Consider the $p\times(q-1)$ matrix $B$ having columns $a_1,\ldots,a_{q-1}$.
Then, the matrix $A^\trp A$ can be written in partitioned form as
\[
A^\trp A\,=\, \left[\begin{array}{cc} B^\trp B & B^\trp a_q\\ a_q^\trp B & a_q^\trp a_q\end{array}\right].
\]
So, by a well-known formula for the determinant of a partitioned positive definite matrix,
\[
\det\bigl(A^\trp A\bigr)\,=\,\det\bigl(B^\trp B\bigr)\cdot\Bigl(a_q^\trp a_q - a_q^\trp B(B^\trp B)^{-1}B^\trp a_q\Bigr).
\]
The second factor on the r.h.s. of the latter equation is equal to \ ${\rm dist}^2(a_q,V_{q-1})$. We have thus obtained that
\ $\det\bigl(A^\trp A\bigr)=\det\bigl(B^\trp B)\cdot{\rm dist}^2(a_q,V_{q-1})$. Now the asserted formula follows by induction on $q$.
\end{itemize}

\end{appendix}



\end{document}